%% file: main.tex
\documentclass[10pt]{article}
\usepackage{amsmath}
\usepackage{graphicx}
\usepackage[inline]{enumitem}
\usepackage{natbib}
\usepackage{url} % not crucial - just used below for the URL
\usepackage{bbm}
\usepackage{bm}
\usepackage{amsfonts}
\usepackage{amsthm}
\usepackage{authblk}
\usepackage{hyperref}
\usepackage[verbose=true]{geometry}

\newgeometry{
    textheight=9in,
    textwidth=5.5in,
    top=1in,
    headheight=12pt,
    headsep=25pt,
    footskip=30pt
}

\newcommand{\E}{E}
\newcommand{\Ehat}{\hat{E}}
\newcommand{\Cov}{\mathrm{cov}}
\newcommand{\Var}{\mathrm{var}}
\newcommand{\indep}{\perp\!\!\!\perp}
\newcommand{\n}[2]{\mathcal{N}\left( #1,#2\right)}
\graphicspath{{plots/}}

\DeclareMathOperator*{\argmin}{arg\,min}

\newtheorem{theorem}{Theorem}%[section]
\newtheorem{corollary}{Corollary}[theorem]
\newtheorem{lemma}{Lemma}[theorem]

\newtheorem{remark}{Remark}
\newtheorem{example}{Example}

\newcommand{\mytitle}{Parameterising the effect of a continuous treatment using average derivative effects}
\title{\mytitle}
\newcommand{\myappendix}{Appendix } 

\author[1]{Oliver J.~Hines}
\author[2]{Karla Diaz-Ordaz}
\author[3]{Stijn Vansteelandt}
\affil[1]{Columbia University, New York, NY, USA}
\affil[2]{University College London, London, U.K.}
\affil[3]{Ghent University, Ghent, Belgium}

\begin{document}
\maketitle

\abstract{
The average treatment effect (ATE) is commonly used to quantify the main effect of a binary treatment on an outcome. Extensions to continuous treatments are usually based on the dose-response curve or shift interventions, but both require strong overlap conditions and the resulting curves may be difficult to summarise. We focus instead on average derivative effects (ADEs) that are scalar estimands related to infinitesimal shift interventions requiring only local overlap assumptions. ADEs, however, are rarely used in practice because their estimation usually requires estimating conditional density functions. By characterising the Riesz representers of weighted ADEs, we propose a new class of estimands that provides a unified view of weighted ADEs/ATEs when the treatment is continuous/binary. We derive the estimand in our class that minimises the nonparametric efficiency bound, thereby extending optimal weighting results from the binary treatment literature to the continuous setting. 
We develop efficient estimators for two weighted ADEs that avoid density estimation and are amenable to modern machine learning methods, which we evaluate in simulations and an applied analysis of Warfarin dosage effects.
}

\input{manuscript}
\bibliography{refs}
\bibliographystyle{apalike}
\pagebreak
\appendix
\input{supplement}
\end{document}

%% file: manuscript.tex
\section{Introduction}

One of the central goals of Biostatistics is to establish the main effect of a treatment/exposure $A$ on an outcome $Y$ and to determine its magnitude and direction. In causal inference, estimands are often defined with reference to potential outcomes $Y^a$ that denote the outcome that would be observed if an intervention had assigned treatment $A=a$. When $A $ is binary, i.e. $A\in \{0,1\}$, the average treatment effect (ATE) $\E(Y^1 - Y^0)$ is a canonical choice \citep{Rosenbaum1983}, but when $A$ is continuous, e.g. dose, duration, frequency, there are several competing extensions. We focus on one such extension based on the average derivative effect (ADE), previously studied in the causal literature by \cite{Rothenhausler2019} and \cite{Athey2021}. We motivate ADEs by arguing that they overcome overlap issues related to alternatives, such as dose-response curves and average shift interventions, and by demonstrating that weighted ADEs/ATEs can be viewed under a common framework based on Riesz's representation theorem.

The dose-response curve $a \mapsto \E(Y^a)$ considers the effect of uniform interventions on a population \citep{Robins2001,Neugebauer2007,Kennedy2017,Hudson2023}, however, inference for such interventions is only possible when all treatment units have a nonzero probability density of receiving treatment $A=a$. Violations of this overlap condition may occur in practice, for example, if an adult drug dose is unsafe to administer to pediatric patients, in which case the effect of intervening at adult dose levels is scientifically uninteresting and impossible to estimate without substantial extrapolation. Such violations can be avoided by examining the dose-response curve in subpopulations where overlap is expected to hold, e.g. pediatric vs. adult patients, but this is rarely done systematically and results in multiple dose-response curves that may be difficult to interpret. 

Similar, albeit weaker concerns hold for the curve of average shift interventions $\varepsilon \mapsto \E(Y^{A+\varepsilon} - Y^A)$, which quantifies the change in mean outcome when the natural value of treatment is shifted by $\varepsilon$ \citep{Munoz2012}. Specifically, when $|\varepsilon|$ is large, the shifted treatment $A+\varepsilon$ may not be supported by the data, and therefore the average effect of a shift intervention is not identified. This concern is alleviated by imposing a threshold on the intervened value of treatment, at the expense of considering more complicated interventions \citep{diaz_nonparametric_2023}. Moreover, dose-response and average shift interventions curves are both univariate functions, with no clear way to summarize the resulting function once it has been obtained

Instead, we consider causal ADE estimands, such as $\lim_{\varepsilon \to 0} \E(Y^{A+\varepsilon} - Y^A)/\varepsilon$ that quantify the effect of shifting treatment by an infinitesimal amount around the natural value of treatment, thus indicating whether increasing/decreasing treatment is generally beneficial. By focusing on realistic treatment values for each individual, ADEs can be identified under weak local overlap assumptions, e.g. when the conditional density of treatment given covariates $X$ is almost surely continuous in treatment. ADEs are also scalar estimands, with no additional summarization required. Under identification assumptions described in Section \ref{sect:prelims}, the ADE and ATE are identified by $\E\{\mu^\prime(A, X)\}$ and $\E\{\mu(1, X) - \mu(0, X)\}$, where $\mu(a, x) \equiv \E(Y|A=a, X=x)$, and superscript prime denotes the derivative w.r.t the first argument.

In modern presentations, such estimands are often viewed as bounded linear functionals $\theta(\mu)$ of an unknown conditional response function $\mu$, see e.g. \cite{hirshberg_augmented_2021} and \cite{chernozhukov_automatic_2024}. A fundamental property used to characterise these is that, for a unique estimand-specific function $\alpha$, one can write $\theta(\mu) = \E\{\alpha(A, X) \mu(A, X)\} = \E\{\alpha(A, X) Y\}$. The representing function $\alpha$, also called the Riesz representer, plays a central role in nonparametric inference for $\theta(\mu)$, with significant recent attention towards directly estimating $\alpha$ from data. Efforts have focused on two estimand subclasses: (a) where $\theta(\mu)$ is the mean of a known linear functional of $\mu$; and (b) where $\alpha(A, X) > 0$ (almost surely) and $\E\{\alpha(A, X)\} = 1$, thus $\alpha$ is a density ratio, also called an importance weight.
For class (a), empirical risk minimization approaches were proposed by \cite{chernozhukov_automatic_2024}, building on the balancing methods of \cite{hainmueller_entropy_2012}, \cite{imai_covariate_2014}, and \cite{zubizarreta_stable_2015}. For class (b), density ratio learning is an active area of research in computer science, with early proposals by \cite{sugiyama_direct_2007} and \cite{kanamori_least-squares_2009}, and a recent review by \cite{hines_learning_2025}.

Although classes (a) and (b) cover a broad variety of scientifically interesting estimands, there remain bounded linear functionals of substantial practical interest that are not contained in either class. For instance, the ATE and the ADE are examples of class (a), but ATEs/ADEs with data-adaptive weights, such as the propensity overlap-weighted treatment effects of \cite{Crump2009} and the density weighted ADEs of \cite{Powell1989} are neither examples of (a) nor (b).

One of the main contributions of the current paper, Section \ref{sect:riesz}, is to introduce and study a new class of bounded linear functional estimands (c) where $\E\{\alpha(A, X) A\} = 1$ and $\E\{\alpha(A, X)\mid X\} = 0$ (almost surely). To give some intuition for this class, consider that, when $\mu(a, x) = \beta a + \omega(x)$ is partially linear, the restrictions of class (c) imply that $\E\{\alpha(A, X) Y \} = \beta$ recovers the linear coefficient of treatment regardless of whether $A$ is binary/continuous. In the general case where $\mu$ is unrestricted, we show that the proposed class corresponds to weighted ATEs when $A$ is binary, and weighted ADEs when $A$ is continuous.
Thus, estimands in class (c) admit causal interpretations as weighted ATEs/ADEs. Although these weighted effects have a more nuanced interpretation than their unweighted counterparts, the latter can impose substantially stronger rate conditions due to the need for density weighting. In settings with high-dimensional confounding, some of the weighted ATEs/ADEs we consider therefore represent more realistic and practically identifiable targets.
When the weights are nonnegative and normalized, the resulting estimands provide interpretable information about the magnitude and direction of the effect of $A$ on $Y$. This is particularly useful for exploratory analyses where no specific intervention is planned: the effect pertains to the corresponding weighted study population, whose characteristics can be described explicitly through weighted summary statistics. More generally, weighted ATEs/ADEs can be used to test the causal null that $Y^a$ does not depend on $a$, since $\theta(\mu) = 0$ under the null.

Furthermore, since class (c) provides a unified view of weighted ATEs/ADEs, we are able to extend the optimal ATE weighting results of \cite{Crump2006, Crump2009} to new optimal ADE weighting results in the setting where $A$ is continuous, with the main additional subtlety being that ADE weights also depend on the treatment, Section \ref{crump_type_result}. More concretely, we derive the estimand in our class that is optimally efficient, in the sense of minimising the nonparametric efficiency bound of an efficient estimator with respect to a sample analogue of $\theta(\mu)$, which is the same definition of efficiency considered by \cite{Crump2009}. When $A$ is binary our optimal estimand reduces to theirs, but when $A$ is continuous our estimand is a new optimally efficient weighted ADE. Although optimally weighted estimands might be criticised for `moving the goalposts', we contend that such estimands offer a pragmatic alternative in low-power settings, e.g. with small sample sizes, low signal to noise ratios, high-dimensional confounding or large imbalance; and in studies based on convenience samples whose relation to any well-defined population is unclear.

Finally, in the second half of our paper, Section \ref{sect_estimation}, we focus on debiased estimators for two estimands of our class, which we call `least squares estimands', due to their connections with nonparametric model projections \citep{Neugebauer2007,Chambaz2012}. These estimands are
\begin{align}
\psi &= \E\left\{ \frac{\Cov(A,Y|X)}{\Var(A|X)} \right\} , \label{lam}  \\
\Psi &=  \frac{\E\left\{ \Cov(A,Y|X) \right\}}{\E\left\{ \Var(A|X) \right\}}.  \label{lam_bar}
\end{align}
Both estimands are connected to our optimally efficient estimand under simplifying assumptions that aid interpretability.
We compare asymptotically efficient, nonparametric estimators of $\psi$ and $\Psi$, the former being a contribution of our work, and the latter following from existing results. These estimators do not require estimation of the conditional treatment density, thus alleviating concerns regarding density learning and kernel estimation of other weighted ADEs, Section \ref{sect:related}. Our estimators are amenable to data adaptive / machine learning of requisite statistical functionals, as we demonstrate on simulated data in Section \ref{sect_simstud}, and on clinical data to determine the effect of Warfarin dose on blood clotting function in Section \ref{sect_appliedIll}.

\section{Related work}
\label{sect:related}

By analyzing weighted ADEs/ATEs nonparametrically, our work complements theoretical analyses of semi-parametric regression estimators in terms of `implied weights'. Specifically, \cite{aronow_does_2016} show that regression estimators induce a reweighting of the original data sample, meaning that, e.g. ATE estimators based on linear models should be viewed as averages of the conditional ATE over reweighted `effective' samples. \cite{chattopadhyay_implied_2023} extend these results by deriving closed form individual weights for several regression estimators, and \cite{Buja2019} show that ordinary least squares coefficients can be interpreted as weighted sums of ‘slopes’ between pairs of observations, without invoking differentiability. In our work, we take a fundamentally different approach by considering population weighted estimands without reference to specific models or estimators, thereby providing control over the choice of weights.

ADEs are also called average partial effects and were originally motivated by parameter estimation in index models in economics \citep{Hardle1989,Powell1989,Newey1993, Imbens2009}. ADE estimation is challenging even when a parametric model for $\mu$ is assumed, since, for all but the simplest parametric models, the estimand is a nonlinear function of the parameters and additional debiasing is required \citep{Hirshberg2018, Wooldridge2020}. Debiased and efficient estimators were proposed by \cite{Newey1993}, using the Riesz representer of the weighted ADE when the weights are known \citep{Powell1989}. These estimators typically use kernel methods to estimate the density of $A$ given $X$ \citep{Hardle1989,Newey1994_mcfaden,Cattaneo2010}. 
However, kernel-based estimators are sensitive to the choice of bandwidth and the asymptotic linearity of the estimator breaks down when the bandwidth is too small \citep{Cattaneo2013}.
Previous studies have avoided density estimation issues by approximating Riesz representers directly, e.g. using the augmented minimax procedure of \cite{hirshberg_augmented_2021}, or the `Riesz regression' procedure of \cite{chernozhukov_automatic_2024} and its variants \citep{hines_automatic_2025,hines_learning_2025}. We instead choose weights that facilitate inference, demonstrating that $\Psi$ and $\psi$ in \eqref{lam} and \eqref{lam_bar} are weighted ADEs that do not require density estimation.

The estimands $\psi$ and $\Psi$ have both appeared in previous literature: $\Psi$ appears in the study of partially linear models \citep{Vansteelandt2020,Chernozhukov2018, Newey2018}; its numerator is used for conditional independence testing \citep{Shah2018}; $\psi$ has been used to estimate the ADE under conditionally linear model assumptions, assuming the semi-parametric model holds \citep{Hirshberg2018}; when $A$ is binary and $X$ is sufficient to adjust for confounding, $\psi$ and $\Psi$ respectively identify the ATE and the propensity overlap weighted effect of treatment on outcome \citep{Crump2006,Crump2009,Robins2008,Li2018,Kallus2020}. In our work, we estimate $\psi$ and $\Psi$ without model assumptions, the former yielding novel estimators.

\section{Weighted average derivative effects}
\label{sect:wade}

\subsection{Preliminaries}
\label{sect:prelims}

Here we introduce weighted ADEs, also called incremental treatment effects by \cite{Rothenhausler2019}, as causal estimands that quantify the effect of an incremental shift in a continuous treatment on an outcome.
Consider $n$ iid observations, $(o_1,...,o_n)$ of a random variable $O=(Y,A,X)$ distributed according to an unknown distribution $P_0$, where $Y\in \mathbb{R}$ is an outcome, $A\in \mathbb{R}$ is a continuous `treatment' of interest and $X \in \mathbb{R}^p$ is a $p$-dimensional vector of covariates. Define: $\mu(A,X)\equiv\E(Y|A,X)$; $\mu(X) \equiv \E(Y|X)$; $\sigma^2(A,X)\equiv \Var(Y|A,X)$; $\pi(X) \equiv \E(A|X)$; $\beta(X)\equiv \Var(A|X)$; $\lambda(X) \equiv \Cov(A, Y|X) / \beta(X)$; $p_0(a|x)$ is the conditional density of $A$ given $X$; $p_0(a,x)$ is the joint density of $(A, X)$; $p_0(a)$ is the marginal density of $A$; and let superscript prime denote the derivative with respect to the first argument, e.g. $\mu^\prime(a, x)$ is the derivative of $\mu$ w.r.t. $a$. Also, let $Y^a$ denote the outcome that would be observed if treatment had taken the value $A=a$. We define the relative conditional shift effect at $a$ for a shift $\varepsilon \neq 0$ as
\begin{align*}
    \delta_{\varepsilon}(a,x) \equiv \frac{\E(Y^{a+\varepsilon} - Y^a | X=x)}{\varepsilon}
\end{align*}
which represents the difference in the conditional mean outcome (per unit treatment) in two counterfactual worlds, where all treatment units receive treatment $a+\varepsilon$ and $a$ respectively. We define the weighted average shift effect as $\Theta_{w,\varepsilon} \equiv \E\{w(A,X)\delta_{\varepsilon}(A,X)\}$, where $w: \mathbb{R} \times \mathbb{R}^p \to \mathbb{R}$ is a weight function such that $k\equiv\E\{w(A,X)\}$ is finite and non-zero. We say that the weight is `normalised' when $k=1$, and for the purposes of interpretation it is often desirable for the weights to be non-negative, though for full generality, we do not impose this as a restriction. When $w(A,X) = 1$, the average shift effect reduces to $\Theta_{1,\varepsilon} = \E(Y^{A+\varepsilon} - Y^A) / \varepsilon$.

Denoting the limit $\delta_{0}(a,x) \equiv \lim_{\varepsilon\to 0} \delta_{\varepsilon}(a,x)$, which is assumed to exist, the weighted ADE is $\Theta_{w,0} \equiv \E\{w(A,X)\delta_{0}(A,X)\}$ and the unweighted ADE refers to $\Theta_{1,0}$. The weighted ADE therefore summarises the curve $\varepsilon \mapsto \Theta_{w,\varepsilon}$ and overcomes the key overlap issue that, when $|\varepsilon|$ is large, $A+\varepsilon$ may be outside of the support of the treatment distribution. Hence, ADEs are appealing summaries of the weighted average shift effect that focuses on realistic (modest) treatment interventions and quantifies their effect independently of the shift size. Moreover, $\Theta_{w,0} = 0$ under the causal null that $Y^a$ does not depend on $a$, regardless of the choice of weight.

To identify the weighted ADE we assume: (i) $Y^a\indep A \mid X$; (ii) $A=a \implies Y^a=Y$; (iii) $p_0(A|X)$ is continuous in $A$ almost surely; (iv) $\mu^\prime(a, x)$ exists. Under these assumptions $\Theta_{w,0} = \theta_{w}(\mu) \equiv \E\{w(A,X)\mu^\prime(A, X)\}$. A proof of this result is provided in \myappendix \ref{proof:identification}, which follows from similar results by \citet[Proposition 1]{Rothenhausler2019}.

\begin{remark}
    Assumption (iii) is a local overlap assumption, since it requires that the treatment density is nonzero for treatment values in the neighborhood of the observed treatment $A$. This is much weaker than the analogous overlap assumption required to identify the dose-response curve: that $p_0(a|X) > 0$ almost surely for each $a$ on the dose-response curve.
\end{remark}

\begin{remark}
Assumptions (i) and (ii) are also used to identify the weighted ATE with (iii) replaced with the binary analogue, $P(A=a|X) > 0$ almost surely, for $a\in\{0,1\}$. The analogous identification result is $\E[w(X)\{Y^1-Y^0\}] = \E[w(X)\{\mu(1,X) - \mu(0,X)\}]$.
\end{remark}

To further develop the theory for weighted ADEs, we consider a Hilbert space $\mathcal{H}$ of functions $f:\mathbb{R} \times \mathbb{R}^{p} \to \mathbb{R}$ equipped with inner-product $\langle f , g \rangle \equiv \E\{f(A,X)g(A,X)\}$; norm $\|f\| = \langle f, f \rangle^{1/2}$; and we assume that $\mu \in \mathcal{H}$, i.e. that $\|\mu\| < \infty$. The Riesz representation theorem implies that, for a bounded continuous linear functional $\theta: \mathcal{H} \to \mathbb{R}$, there exists a unique function $\alpha \in \mathcal{H}$, called the Riesz representer, such that $\theta(f) = \langle f, \alpha\rangle$ for all $f \in \mathcal{H}$. Hence, $\theta(\mu) = \langle \mu, \alpha\rangle = \E\{\alpha(A, X) Y\}$.

The weighted ADE $f \mapsto \theta_w(f)$ is a bounded linear functional under limited regularity conditions on $P_0$ and $w$. Specifically, we assume: (C1) the support of $A$ given $X=x$, denoted $\Omega(x)$, is an open (possibly unbounded) interval; (C2) $w(a, x)p_0(a|x)$ is differentiable w.r.t. $a$ for all $a \in \Omega(x)$; (C3) $\alpha_w \in \mathcal{H}$, and hence the norm $\|\alpha_w\| < \infty$ is finite, where
\begin{align}
\alpha_w(a,x) &= -w^\prime(a,x) - w(a,x) \frac{p_0^\prime(a|x)}{p_0(a|x)} \label{canonical_contrast}
\end{align}
is the Riesz representer of $\theta_w(f)$. Under these assumptions $\theta_w(f) = \langle f, \alpha_w\rangle$ is an inner product on $\mathcal{H}$, and hence $\theta_w$ is a bounded linear functional with Riesz representer $\alpha_w$. This result follows from integration by parts, see \myappendix \ref{append:integration_by_parts} and \cite{Powell1989} for details. Thus, for the weighted ADE, $\theta_w(\mu) = \langle \mu, \alpha_w \rangle = \E\{\alpha_w(A,X) Y\}$. In the next section, we study weighted ADEs using a set of functions $\mathcal{R} \subseteq \mathcal{H}$ such that $\langle \mu, \alpha \rangle$ is a weighted ADE for all $\alpha \in \mathcal{R}$. This essentially involves treating \eqref{canonical_contrast} as an ordinary differential equation and solving for $w$.

\begin{remark}
    \label{remark:prelim}
    In the setting where $A$ is a binary treatment, the analogous Riesz representation of the weighted ATE is $\E[w(X)\{\mu(1,X)-\mu(0,X)\}] = \langle \mu, \alpha \rangle$, where $\alpha(a,x) = w(x)\{a-\pi(x)\}/ [\pi(x)\{1-\pi(x)\}]$ is a Riesz representer. Weighted ATEs are said to be normalised when $\E\{w(X)\}=1$.
\end{remark}

\subsection{A Riesz representer class}
\label{sect:riesz}

Consider the estimands $\theta(\mu) = \langle \mu, \alpha \rangle$ where $\alpha$ is an element of the set 
\begin{align*}
    \mathcal{R} \equiv \left\{\alpha \in \mathcal{H} ~ \Big| ~
        \E\{\alpha(A,X)A\} = 1 ,~
        \E\{\alpha(A,X)|X\} = 0
    \right\}.
\end{align*}
This class of estimands is motivated by the fact that, when $k = 1$, the Riesz representer $\alpha_w$ in \eqref{canonical_contrast} is a member of $\mathcal{R}$, i.e. normalised weighted ADEs belong to our class of estimands. Moreover, when $A$ is a binary treatment, then the Riesz representer of the normalised weighted ATEs (Remark \ref{remark:prelim}) is also a member of $\mathcal{R}$, suggesting that a unified understanding of weighted ATEs and weighted ADEs may be obtained by studying $\mathcal{R}$.

Through Theorem \ref{alpha_theorem} we show that, for each $\alpha \in \mathcal{R}$, one can construct a weight $w$ such that $\theta_w(\mu) = \langle \mu, \alpha \rangle$ is a weighted ADE. This implies an isomorphism between weighted ADEs and estimands of our class. There is no guarantee, however, that the weight in Theorem \ref{alpha_theorem} is non-negative. We address this by deriving a sufficiency condition for weight non-negativity (Lemma \ref{sufficiency}). The significance of these results is that, under suitable identification assumptions, any estimand of the form $\langle \mu, \alpha \rangle$ can be ascribed a causal interpretation in terms of weighted ADEs. Moreover, new weighted ADEs can be specified by their Riesz representers rather than their weight functions, a fact that we exploit with reference to the optimal weighting strategies in Section \ref{crump_type_result}. 

\begin{theorem}
\label{alpha_theorem}
Let $F(a|x)$ be the distribution function of $A$ given $X=x$ and assume that the support of $A$ given $X=x$ is an open (possibly unbounded) interval. For $\alpha \in \mathcal{R}$ define the weight
\begin{align}
 w(a,x) = \frac{F(a|x)\{1-F(a|x)\}}{p_0(a|x)}\left[\E\{\alpha(A,X)|A > a, X=x\} - \E\{\alpha(A,X)|A \leq a, X=x\}\right]. \label{f_tilde_theorem_alpha}
\end{align}
For all differentiable functions $f \in \mathcal{H}$, $\langle f, \alpha \rangle = \E\{w(A,X)f^\prime(A,X)\}$.
Proof in \myappendix \ref{main_theorem_proof}.
\end{theorem}

\begin{lemma}
\label{sufficiency}
If $\alpha(a,x)$ is monotonically increasing in $a$ then the weight in \eqref{f_tilde_theorem_alpha} is non-negative.\\
For proof note that $\E\{\alpha(A,X)|A > a, X=x\} \geq \alpha(a, x) \geq \E\{\alpha(A,X)|A \leq a, X=x\}$ under monotonicity, and hence each component of \eqref{f_tilde_theorem_alpha} is non-negative.
\end{lemma}

\begin{remark}
\label{alpha_remark}
In the setting where $A$ is a binary treatment, then an analogue of Theorem \eqref{alpha_theorem} is obtained by letting $w(x) = \E\{\alpha(A,X)A|X=x\}$, with $\langle f, \alpha \rangle = \E[w(X)\{f(1,X) - f(0,X)\}]$. Thus, any estimand of the form $\langle \mu, \alpha \rangle$ can be ascribed a causal interpretation in terms of weighted ATEs when the treatment is binary.
\end{remark}

We apply Theorem \ref{alpha_theorem} and Lemma \ref{sufficiency} in the following examples, which illustrate the connection between the Riesz representer and the ADE weight. In each example, we assume that the resulting Riesz representer, $\alpha \in \mathcal{H}$, and hence the norm $\| \alpha  \| < \infty$ is finite. Each of these examples identifies a causal effect of treatment on outcome, with weighting determined by the ADE weight.

\begin{example}[Average derivative effect (ADE)]
\label{ade}
Setting $w(a,x) = 1$ results in the ADE, $\E\{\mu^\prime(A,X)\}$ with Riesz representer $\alpha(a,x) = -p_0^\prime(a|x) / p_0(a|x)$, which follows from \eqref{canonical_contrast}. The ADE is normalised since $\E\{\alpha(A,X)A\} = \E\{w(A,X)\}=1$ and hence $\alpha \in \mathcal{R}$.
\end{example}

\begin{example}[Density weighted ADE]
\label{dens_ade}
Density weights are designed to facilitate inference by avoiding the inverse density weighting in the Riesz representer in Example 1 \citep{Powell1989,Cattaneo2010}.
The density weight sets \mbox{$w(a,x)=p_0(a,x)$} with Riesz representer $- 2 p_0^\prime(a,x)$. This estimand is unnormalised with $k=\E\{p_0(A,X)\}$. A normalised density weighted ADE is obtained by letting \mbox{$w(a,x)= p_0(a,x)/\E\{p_0(A,X)\}$} in which case the Riesz representer becomes $\alpha(a,x) = - 2 p_0^\prime(a,x) / \E\{p_0(A,X)\}$ with $\alpha \in \mathcal{R}$.
\end{example}

\begin{example}[Average dose-response derivative]
\label{adrd}
Under standard assumptions, the dose-response curve is identified by $a \mapsto \varphi(a) \equiv\E\{\mu(a,X)\}$. The mean derivative, $\E\{\varphi^\prime(A)\}$ is a weighted ADE with weight $w(a,x)=p_0(a)/p_0(a|x)$, where we make the overlap assumption that $p_0(a|X) > 0 $ almost surely, for all $a$ on the support of $A$, i.e. when $p_0(a) > 0$. By \eqref{canonical_contrast}, the corresponding Riesz representer is $\alpha(a,x) = - p_0^\prime(a) / p_0(a|x)$, with $\E\{\alpha(A,X)A\} = \E\{w(A,X)\}=1$, hence $\alpha \in \mathcal{R}$.
\end{example}

\begin{example}[Least Squares Estimands]
\label{alse}
The estimands $\psi$ and $\Psi$ in \eqref{lam} and \eqref{lam_bar} are both of the form $\langle \mu, \alpha \rangle$ with Riesz representers $\alpha_\psi(a, x) \equiv \{a - \pi(x)\} / \beta(x)$ and $\alpha_\Psi(a, x) \equiv \alpha_\psi(a, x) \beta(x) / \E\{\beta(X) \}$
with $\alpha_\psi, \alpha_\Psi \in \mathcal{R}$. Theorem \ref{alpha_theorem} implies that $\psi$ and $\Psi$ are ADEs with weights
\begin{align}
w_\psi(a,x) &= \frac{F(a|x)\{1-F(a|x)\}}{p_0(a|x)\beta(x)} \{E(A|A>a,X=x) - E(A|A\leq a,X=x) \}  \label{ALSE-weight-new}
\end{align}
and $w_\Psi(a,x) = w_\psi(a,x)\beta(x) / \E\{\beta(X) \}$, respectively. Both weights are non-negative by Lemma \ref{sufficiency}. We call $\psi$ and $\Psi$ least squares estimands due to their connections with the model projections (\myappendix \ref{sect:projection}). In \myappendix \ref{append_ls_weight} we examine in detail how $w_\psi(a,x)$ looks for various parametric treatment distributions. Notably, $w_\psi(a,x) = 1$ when $A$ is normally distributed given $X$, thus, in this case, $\psi$ is the unweighted ADE and $\Psi$ is the unweighted ADE when $\beta(x)$ is constant. Although there is no need to characterise or estimate the treatment weight to use and interpret $\psi$ and $\Psi$ in practice, an approach for estimating $w_\psi$ is outlined in \myappendix \ref{weight_approximations} based on a location-scale model used by \cite{Kennedy2017} and \cite{Klyne2023}.
\end{example}

\begin{example}[Dichotomised treatment]
\label{dco_estimand}
The causal effect of a continuous treatment on an outcome can be quantified via the ATE of the dichotomised treatment $\mathbb{I}(A > a_0)$ which, for a predetermined constant $a_0$, takes the value 1 when $A > a_0$ and 0 otherwise. This approach is common and has previously been discouraged on the grounds of bias, efficiency, and because it is not clear what hypothetical intervention is being considered \citep{Vanderweele2011, Berzuini2013}. However, recent studies have justified the use of dichotomised treatments by characterising them in terms of threshold interventions \citep{van_der_laan_nonparametric_2023} and modified treatment policies \citep{Lee2024}. Under standard causal assumptions for binary treatments, the resulting ATE is identified by
$\E\{\E(Y|A>a_0, X) - \E(Y|A\leq a_0, X)\}/c$,
for a normalisation constant $c$. This estimand is of the form $\langle \mu, \alpha \rangle$ with Riesz representer and implied ADE weight
\begin{align*}
  \alpha(a,x) &= \frac{\mathbb{I}(a > a_0) - \{1 - F(a_0|x)\} }{F(a_0|x)\{1 - F(a_0|x)\} c}, \\
  w(a,x) &= \left\{ \frac{F(a|x)\{1 - F(a|x)\}}{F(a_0|x)\{1 - F(a_0|x)\}}\right\} \frac{P(A>a_0|A>a,X=x) - P(A>a_0|A\leq a,X=x)}{p_0(a|x) c}.
\end{align*}
When $c = \E\{\E(A|A>a_0, X) - \E(A|A\leq a_0, X)\}$ then $\E\{\alpha(A,X)A\} = 1$ and hence $\alpha \in \mathcal{R}$. Moreover, $w(a,x)$ is non-negative by Lemma \ref{sufficiency}.
\end{example}

An interpretable subclass of estimands are those of the form $\langle\mu, \alpha \rangle$ where $\alpha$ is an element of 
\begin{align*}
    \mathcal{R}^* \equiv \left\{\alpha \in \mathcal{H} ~ \Big| ~
        \E\{\alpha(A,X)A|X\} = 1 ,~
        \E\{\alpha(A,X)|X\} = 0
    \right\},
\end{align*}
with $\mathcal{R}^* \subseteq \mathcal{R}$.
% $\mathcal{R}^* \equiv \{\alpha \in \mathcal{R} \mid \E\{\alpha(A,X)A|X\} = 1 \}$.
When $A$ is binary, this subclass contains only the unweighted ATE. However, when $A$ is continuous, it corresponds to weighted ADEs where the weight in Theorem \ref{alpha_theorem} is normalized to $\E\{w(A,X)|X\} = 1$. Examples of weighted ADEs with this property include the unweighted ADE (Example \ref{ade}), and the least squares estimand $\psi$ (Example \ref{alse}). This normalisation aids interpretability, since, within each covariate level $X=x$, the ADE weight $w(a,x)$ provides a normalized reweighting according to the value of treatment, with the total amount of weight given to covariate level $x$ unchanged.

\subsection{Optimally efficient estimands}
\label{crump_type_result}

Having characterised the class of weighted ADE/ATE estimands as $\langle \mu, \alpha \rangle$ for $\alpha \in \mathcal{R}$, this characterisation raises the question: \textit{what estimand in this class has the most favourable efficiency bound?} An analogous study in the context of weighted ATEs was conducted by \cite{Crump2006,Crump2009}, and, in this section, we extend their results to the setting of continuous treatments and weighted ADEs, with the extra subtlety being that ADE weights are functions of the treatment as well as the covariates. In particular, we derive the Riesz representer $\alpha \in \mathcal{R}$ that optimises the efficiency bound of a sample analogue of $\theta_w(\mu)$. This optimality criteria is exactly that of Crump et al., and the connection between their work and ours is made explicit in Remark \ref{crump_remark}. 
% By Theorem \ref{alpha_theorem}, the optimal Riesz representer implies a corresponding optimally weighted ADE. 

We rely on efficient influence curves to characterize the sensitivity of weighted ADEs to small changes in the data distribution. Influence curves are model-free, mean zero, functionals of the data distribution, derived from the definition of the target estimand. They are useful for constructing efficient estimators and deriving their asymptotic efficiency bounds \citep{Hines2021,Fisher2020}. This efficiency bound is a property of the estimand itself and is given by the variance of its influence curve, which is finite.
When the weight $w(a,x)$ is known and (C1) to (C3) are assumed, the influence curve of $\theta_w = \theta_w(\mu)$ is
\begin{align}
\phi_{\theta,w}(o) =  \alpha_w(a,x) \{y-\mu(a,x)\} + w(a,x)\mu^\prime(a,x) -  \theta_w  \label{newey_ic}
\end{align}
where $\alpha_w(a,x)$ is the Riesz representer in \eqref{canonical_contrast} and $o=(y,a,x)$ \citep{Newey1993}. In all but Example \ref{ade} the weight function is unknown, however, the influence curve above, derived when the weight is known, offers some insight into optimal weight selection. Specifically, we consider the efficiency bound of an efficient estimator, $\hat{\theta}_w$, of the sample analogue of $\theta_w$,
\begin{align*}
\theta_{w,S} &\equiv n^{-1}\sum_{i=1}^n w(a_i,x_i)\mu^\prime(a_i,x_i) \\
\sqrt{n}(\hat{\theta}_w-\theta_{w,S}) & \overset{d}{\to} \n{0}{V} \\
V &\equiv \E\{\alpha_w^2(A,X)\sigma^2(A,X)\}.
\end{align*}
The efficiency bound with respect to $\theta_{w,S}$, rather than $\theta_w$, is chosen so that the final two terms in \eqref{newey_ic} may be disregarded. Not only does this simplify the subsequent analysis, but these terms capture the difference between the weighted ADE conditional on the sample distribution and that of the population as a whole, which depends on the unknown value of $\theta_w$. I.e.,
\begin{align*}
\sqrt{n}(\hat{\theta}_w-\theta_{w}) & \overset{d}{\to} \n{0}{V+U} \\
U &\equiv \E[\{w(A,X)\mu^\prime(A,X) -  \theta_w\}^2]
\end{align*}
thus selecting weights to minimise $V+U$ is conceptually problematic as $\theta_w$ is itself the target estimand \citep{Crump2006,Crump2009}. Theorem \ref{alpha_theorem}, offers constraints on $\alpha_w(a,x)$, under which $\alpha_w$ is the Riesz representer of a weighted ADE. Our goal, therefore, is to find $\alpha_w \in \mathcal{R}$ which minimises $V$. The optimal solution is given by Theorem \ref{optimality_1}. For the reasons of interpretability discussed in Section \ref{sect:riesz}, we also find the Riesz representer $\alpha_w \in \mathcal{R}^*$ which minimises $V$, with the solution given by Theorem \ref{optimality_2}. 
The Riesz representers in Theorems \ref{optimality_1} and \ref{optimality_2} each imply a corresponding ADE weight according to Theorem \ref{alpha_theorem}. Under outcome homoscedasticity, as described in Corollaries \ref{corol:indep} and \ref{corol:indep2}, these estimands reduce to the least squares estimands $\Psi$ and $\psi$ receptively, with non-negative ADE weights according to Lemma \ref{sufficiency}. We state both theorems in terms of the quantities
\begin{align*}
    \tilde{\pi}(x) &\equiv \frac{\E\{A / \sigma^2(A,X) \mid X = x\}}{\E\{1 / \sigma^2(A,X) \mid X = x\}}, \\
    d(a, x) &\equiv \frac{a - \tilde{\pi}(x)}{\sigma^2(a, x) }.
\end{align*}

\begin{theorem}
\label{optimality_1}
Minimizing the efficiency bound $V=n\Var\{\hat{\theta}_w-\theta_{w,S}\}$, subject to the constraints $\E\{\alpha_w(A,X)A\}=1$ and $\E\{\alpha_w(A,X)\mid X\}=0$, has the solution $\alpha_w(a,x) = d(a, x) / \E\{d(A, X)A\}$,
which implies the optimally efficient estimand
\begin{align*}
  \frac{\E\left\{d(A,X) Y\right\}}{\E\left\{d(A,X) A \right\}}.
\end{align*}
Proof in \myappendix \ref{optimality_proof}.
\end{theorem}

\begin{corollary}
\label{corol:indep}
When $Y$ is homoscedastic, i.e. $\sigma^2(a,x)$ is constant, the optimal estimand in Theorem \ref{optimality_1} is $\Psi$. For proof, note that $\tilde{\pi}(x) = \pi(x)$ when $\sigma^2(a,x)$ is constant.
\end{corollary}

\begin{remark}
\label{crump_remark}
In the setting where $A$ is a binary treatment, \cite{Crump2006,Crump2009} consider the weighted ATE $\E[w(X) \{\mu(1,X) - \mu(0,X)\}]$ where $w(x)$ is known. They derive that the efficiency bound of a sample analogue  weighted ATE has the same form as $V$ above, with Riesz representer $\alpha(a,x) = w(x)\{a-\pi(x)\}/[\pi(x)\{1-\pi(x)\}]$. They minimise $V$ over $w(x)$ with the normalisation constraint $\E\{w(X)\} = 1$, and their result is recovered by Theorem \ref{optimality_1} when $A$ is binary, see \myappendix \ref{crump_appendix}. Thus, Theorem \ref{optimality_1} represents a generalisation that covers continuous and binary treatments.
\end{remark}

\begin{theorem}
\label{optimality_2}
Minimizing the efficiency bound $V=n\Var\{\hat{\theta}_w-\theta_{w,S}\}$, subject to the constraints $\E\{\alpha_w(A,X)A \mid X\}=1$ and $\E\{\alpha_w(A,X) \mid X\}=0$, has the solution $\alpha_w(a, x) = d(a, x) / \E\{d(A, X)A\mid X=x\}$,
which implies the estimand
\begin{align*}
  \E\left[ \frac{\E\left\{d(A,X)Y \mid X\right\}}{\E\left\{d(A,X)A \mid X\right\}} \right].
\end{align*}
Proof in \myappendix \ref{optimality_proof2}.
\end{theorem}

\begin{corollary}
\label{corol:indep2}
When $Y$ is homoscedastic conditional on $X$, i.e. $\sigma^2(a,x)=\sigma^2(x)$, the estimand in Theorem \ref{optimality_2} is $\psi$. For proof, note that $\tilde{\pi}(x) = \pi(x)$ when $\sigma^2(a,x)=\sigma^2(x)$. 
\end{corollary}

\begin{remark}
    When $A$ is a binary treatment, the estimand in Theorem \ref{optimality_2} reduces to the unweighted ATE, since the constraints on the Riesz representer uniquely imply this estimand.
\end{remark}

It is insightful to consider which parts of the unknown data distribution $P_0$ contribute to the Riesz representers of the optimal estimands in Theorem \ref{optimality_1} and \ref{optimality_2}.
To make such a comparison, we consider that the Markov factorization of $P_0$ consists of three conditional components: $P_{Y|A,X}$, $P_{A|X}$, and $P_X$. The optimal Riesz representer in Theorem \ref{optimality_1} depends on all three components, with the Riesz representer in Theorem \ref{optimality_2} depending on $P_{Y|A,X}$ and $P_{A|X}$. However, the dependence on $P_{Y|A,X}$ makes these estimands challenging to interpret. For example, in studies where there are multiple outcomes of interest, different weighted populations would be considered for each outcome.
This concern is alleviated under the simplifying assumptions on $P_{Y|A,X}$ in Corollaries \ref{corol:indep} and \ref{corol:indep2}, namely: when $\sigma^2(a, x)$ is constant the estimand in Theorem \ref{optimality_1} reduces to $\Psi$, whose Riesz representer depends on $P_{A|X}$ and $P_X$; when $\sigma^2(a, x) = \sigma^2(x)$ the estimand in Theorem \ref{optimality_2} reduces to $\psi$, whose Riesz representer depends on $P_{A|X}$.
Whereas the optimal estimands lack suitable interpretability, the simpler estimands $\Psi$ and $\psi$ are more interpretable, since they are guaranteed to have non-negative ADE weights that do not depend on the outcome distribution, and both estimands can be interpreted within a least squares framework. We therefore recommend focusing on $\Psi$ and $\psi$ instead of the optimal estimands in Theorems \ref{optimality_1} and \ref{optimality_2}.
A similar recommendation is usually made in the binary treatment setting, where it is common to focus on the overlap weighted treatment effect $\Psi$ rather than the optimal estimand in Theorem \ref{optimality_1} \citep{Crump2006}.

\section{Estimation}
\label{sect_estimation}
\subsection{Efficient estimators}

Here we focus on efficient estimation of $\psi$ and $\Psi$ in \eqref{lam} and \eqref{lam_bar}. Nonparametric efficient estimators for $\Psi$ have previously been studied by \cite{Robinson1988}, \cite{Chernozhukov2018}, and \cite{Vansteelandt2020}, with a recent emphasis on the smoothness and convergence rates of nuisance estimators \citep{Newey2018,Balakrishnan2023}. To our knowledge, nonparametric efficient estimators for $\psi$ have not previously been developed, with \cite{Hirshberg2018} instead considering estimation under a semi-parametric model.
The influence curves of $\psi$ and $\Psi$ respectively are
\begin{align*}
\phi_{\psi}(o) &= \frac{\{a-\pi(x)\}}{\beta(x)}\left[y - \mu(x) - \lambda(x)\{a-\pi(x)\} \right] +  \lambda(x)- \psi \\
\phi_{\Psi}(o) &= \frac{\{a-\pi(x)\}}{\E\{\beta(X)\}}\left[y - \mu(x) - \Psi \{a-\pi(x)\}\right].
\end{align*}
These influence curves may be used to construct efficient estimators by setting an estimate of the sample mean influence curve to zero. For $\psi$, this strategy is equivalent to so-called one-step correction outlined in \myappendix \ref{appendix:asymptotics}. The resulting estimators are
\begin{align*}
\hat{\psi} &= n^{-1}\sum_{i=1}^n  \frac{\{a_i - \hat{\pi}(x_i)\}}{\hat{\beta}(x_i)}[y_i - \hat{\mu}(x_i) - \hat{\lambda}(x_i)\{a_i - \hat{\pi}(x_i)\}] +  \hat{\lambda}(x_i), \\
\hat{\Psi} &= \frac{\sum_{i=1}^n\{a_i - \hat{\pi}(x_i)\}\{y_i - \hat{\mu}(x_i)\}}{\sum_{i=1}^n\{a_i - \hat{\pi}(x_i)\}^2},
\end{align*}
where superscript hat denotes consistent function estimators. In practice, we recommend cross-fitting, as described in Section \ref{proposed_algos}, to obtain the fitted models and evaluate the estimators using a single sample \citep{Chernozhukov2018,van_der_laan_cross-validated_2011}. We discuss the reasons for cross-fitting with reference to Theorems \ref{asym_theorem_alse1} and \ref{asym_theorem_alse2}, which give conditions under which $\hat{\psi}$ and $\hat{\Psi}$ are regular asymptotically linear (RAL).

\begin{theorem}
\label{asym_theorem_alse1}
Assume that there exists constants $\epsilon>0, K>0$ such that (almost surely) $\epsilon < \hat{\beta}(X)$, $\beta(X) \in (\epsilon, K) $, $\hat{\lambda}(X)\in(-K,K)$, $\Var(Y|X) < K$, $E[\{Y-\mu(X)\}^4 |X] < K$, $E[\{A-\pi(X)\}^4 |X] < K$. Suppose also that at least one of the following two conditions hold:
\begin{enumerate}[label=\arabic*.]
\item (Sample-splitting) $\hat{\pi}(x), \hat{\mu}(x), \hat{\lambda}(x)$, and $\hat{\beta}(x)$ are obtained from a sample independent of the one used to construct $\hat{\psi}$.
\item (Donsker condition) The quantities $\hat{\lambda}(X)$,
\begin{align*}
\frac{\{A-\hat{\pi}(X)\}\{Y - \hat{\mu}(X)\}}{\hat{\beta}(X)}, \frac{\hat{\lambda}(X)\{A-\hat{\pi}(X)\}^2}{\hat{\beta}(X)} 
\end{align*}
fall within a $P$-Donsker class with probability approaching $1$.
\end{enumerate}
Finally, letting $\|.\|$ denote the $\mathcal{L}_2(X)$ norm, assume
\begin{enumerate}[label=(A\arabic*)]
\item $\|\pi-\hat{\pi}\| = o_P(n^{-\nu/4})$ and $\|\mu-\hat{\mu}\| = o_P(n^{-\tau/4})$ where $\nu \geq1$, $\tau\geq0$ and $\nu + \tau \geq 2$.
\item The product of $\|\lambda-\hat{\lambda}\|$ and $\|\beta-\hat{\beta}\|$ is $o_P(n^{-1/2})$.
\end{enumerate}
Then $\hat{\psi}$ is RAL with influence curve, $\phi_{\psi}(O)$, and hence $\sqrt{n}(\hat{\psi}-\psi)$ converges in distribution to a mean-zero normal random variable with variance $\E\{\phi_{\psi}^2(O)\}$.
Proof in \myappendix \ref{appendix:asymptotics}.
\end{theorem}

\begin{theorem}
\label{asym_theorem_alse2}
Assume (A1) in Theorem \ref{asym_theorem_alse1}, the quantity $n^{-1}\sum_{i=1}^n\{a_i - \hat{\pi}(x_i)\}^2 > 0$, and $\E[\{A-\pi(X)\}^2] > 0$, there exists a constant $K>0$ such that $\Var(Y|X) < K$ and $\beta(X) < K$ and suppose that at least one of the following two conditions hold:
\begin{enumerate}[label=\arabic*.]
\item (Sample-splitting) $\hat{\pi}(x)$ and $ \hat{\mu}(x)$ are obtained from a sample independent of the one used to construct $\hat{\Psi}$.
\item (Donsker condition) The quantities $\{A-\hat{\pi}(X)\}\{Y - \hat{\mu}(X)\}$ and $\{A-\hat{\pi}(X)\}^2$ fall within a $P$-Donsker class with probability approaching $1$.
\end{enumerate}
Then $\hat{\Psi}$ is RAL with influence curve, $\phi_{\Psi}(O)$, and hence $\sqrt{n}(\hat{\Psi}-\Psi)$ converges in distribution to a mean-zero normal random variable with variance $\E\{\phi_{\Psi}^2(O)\}$.
Proof in \myappendix \ref{appendix:asymptotics}.
\end{theorem}

The estimator $\hat{\psi}$ requires learning $\beta$ and $\lambda$, whereas $\hat{\Psi}$ does not, with Theorem \ref{asym_theorem_alse1} requiring (A2) to control the error in $\hat{\beta}$ and $\hat{\lambda}$. This distinction makes $\Psi$ generally more straightforward to estimate than $\psi$. Assumption (A2) also demonstrates that $\hat{\psi}$ is `rate double robust', in the sense that $\hat{\lambda}$ can converge slowly, so long as $\hat{\beta}$, converges sufficiently quickly, and vice-versa. Thus, one can trade-off accuracy in $\hat{\lambda}$ and $\hat{\beta}$. Similar double robustness has been demonstrated previously, e.g. for the augmented inverse probability weighted (AIPW) estimator of the ATE \citep{Robins1994}, which trades-off accuracy in the propensity score and outcome estimators. On top of this double robustness, (A1) implies a one-sided robustness of $\hat{\psi}$ and $\hat{\Psi}$ with respect to the estimators $\hat{\pi}$ and $\hat{\mu}$. Specifically, $\hat{\mu}$ can converge slowly, so long as $\hat{\pi}$ converges sufficiently quickly, but the converse is not true, since (A1) requires that $\hat{\pi}$ converges at least at $n^{1/4}$ rate.

The Donsker conditions in Theorems \ref{asym_theorem_alse1} and \ref{asym_theorem_alse2} are usually not guaranteed to hold when flexible machine learning methods are used to estimate nuisance functions. Fortunately, cross fitting of nuisance functions offers a way of avoiding Donsker conditions, at the expense of making nuisance functions more computationally expensive to learn \citep{Chernozhukov2018,van_der_laan_cross-validated_2011}. 

\begin{remark}
When $A\in\{0,1\}$ is a binary treatment, $\hat{\psi}$ reduces to the AIPW estimator of the ATE. In particular, if we estimate $\hat{\mu}(a,x)$ and $\hat{\pi}(x)$, then let $\hat{\lambda}(x)=\hat{\mu}(1,x)-\hat{\mu}(0,x)$, $\hat{\beta}(x)=\hat{\pi}(x)\{1-\hat{\pi}(x)\}$, and $\hat{\mu}(x) =\hat{\mu}(0,x) + \hat{\lambda}(x)\hat{\pi}(x)$, one obtains the AIPW estimator
\begin{align*}
n^{-1}\sum_{i=1}^n  \frac{\{a_i - \hat{\pi}(x_i)\}}{\hat{\pi}(x_i)\{1-\hat{\pi}(x_i)\}}\{y_i - \hat{\mu}(a_i,x_i)\} +  \hat{\mu}(1,x_i)-\hat{\mu}(0,x_i).
\end{align*}
Hence, $\hat{\psi}$ is a generalisation of the AIPW estimator that covers binary and continuous treatments.
\end{remark}

\subsection{Nuisance function estimators}

The estimator $\hat{\Psi}$ is indexed by the choice of estimator for $\hat{\mu}$ and $\hat{\pi}$, with the estimator $\hat{\psi}$ additionally indexed by the choice of estimator for $\hat{\lambda}$ and $\hat{\beta}$. Generally, we are not constrained to any particular learning method, making these estimators amenable to data adaptive / machine learning estimation of these working models.
Data adaptive regression algorithms are well developed for the regularised regression of an observed variable on to a set of explanatory variables, e.g. for $\mu$, and $\pi$ in the present context, which can be estimated by respectively regressing $Y$ and $A$ on $X$. For $\lambda$ and $\beta$, however, estimation methods are less well developed, and we propose so-called meta-learning approaches, which estimate $\lambda$ and $\beta$ by solving a series of regression problems.

In the setting where $A$ is a binary treatment and $X$ is sufficient to adjust for confounding, $\lambda$ identifies the conditional ATE, estimation of which is an active area of research, with an emphasis on flexible machine learning methods \citep{Abrevaya2015,Athey2016,Nie2017,Kallus2018,Wager2018,Kunzel2019,Kennedy2020}. Also in the binary treatment setting, $\beta(x) = \pi(x)\{1-\pi(x)\}$ therefore there is no need for a separate estimator of $\beta$. However, the problem of estimating conditional variance functions, more generally, has received some attention in the literature, with applications in constructing confidence intervals for the conditional mean function $\pi$ and for estimating signal-to-noise ratios \citep{Shen2020,Wang2008,Cai2009,Verzelen2018}. We will consider two approaches for estimating $\lambda$ and $\beta$.

The first approach, which we refer to as the direct learning approach, involves decomposing $\lambda$ and $\beta$ into functions of conditional expectations, each of which can be estimated using standard regression methods, with the estimates combined to produce $\hat{\lambda}$ and $\hat{\beta}$. Specifically, letting $\Ehat\{YA|X=x\}$ and $\Ehat\{A^2|X=x\}$ denote estimates obtained by respectively regressing $YA$ and $A^2$ on $X$, we define the direct nuisance estimators
\begin{align}
\hat{\beta}(x)   &= \Ehat\{A^2|X=x\} - \hat{\pi}^2(x). \label{dir_beta} \\
\hat{\lambda}(x) &= \frac{\Ehat\{YA|X=x\} - \hat{\mu}(x)\hat{\pi}(x)}{\hat{\beta}(x)}. \label{dir_lambda}
\end{align}
The issue with this direct approach is that, although regularization methods can be used to control the smoothness of each individual regression function, there is no guarantee on the smoothness of the combined functions $\hat{\lambda}$ and $\hat{\beta}$. In practice, these may be erratic due to artefacts of the regularization of the individual regression functions and there is also no guarantee that $\hat{\beta}$ is greater than zero. Additionally, when $\lambda$ or $\beta$ are simple functions, e.g. constant, one might hope that they would be easy to learn at fast rates of convergence. However, the corresponding direct estimators may inherit slow convergence rates from the estimators $\hat{\mu}, \hat{\pi}, \ldots$, which could be the case when these functions are, in fact, complex. These issues motivate an alternative approach where the smoothness of $\hat{\lambda}$ and $\hat{\beta}$ can be controlled directly, and one can ensure that $\hat{\beta}(x)>0$.

The second approach, which we refer to as quasi-oracle learning, is a meta-learning method based on the R-learner of the conditional ATE $\lambda$ due to \cite{Nie2017}. We begin by characterising a ratio of conditional expectations as a population risk minimiser in the following Lemma, which suggests that such ratios can be estimated by empirical risk minisation.

\begin{lemma}
\label{R_identLemma}
Let $O=(X,V,W)$ be a random variable consisting of $X\in\mathbb{R}^p$, $V\in\mathbb{R}$, and  $W\in\mathbb{R}$ with $W > 0$ almost surely. Let $\mathcal{F}$ denote the set of all possible functions $g:\mathbb{R}^p\to \mathbb{R}$. Then
\begin{align}
\frac{\E(V|X=x)}{\E(W|X=x)} = \argmin_{g \in \mathcal{F}} \E\left[ W\left\{ \frac{V}{W} - g(X) \right\}^2\right]
\label{R-identity}
\end{align}
where the requisite moments of $O$ are assumed to be finite. See \myappendix \ref{appendix:optimal_derivs} for proof.
\end{lemma}

The R-learner is based on Lemma \ref{R_identLemma}, where $W=\{A-\pi(X)\}^2$, $V=\{Y-\mu(X)\}\{A-\pi(X)\}$, and the target ratio is $\lambda(x)$, hence an estimator for $\lambda$ is obtained by regressing $\{Y-\mu(X)\}/\{A-\pi(X)\}$ on $X$ with weights $\{A-\pi(X)\}^2$ using a weighted mean-squared error loss. We call this an `oracle estimator' for $\lambda$, since it is the regression problem that one would like to solve if these `outcomes' and `weights' were observed. The R-learner of the conditional ATE mimics the oracle learner by first estimating $\mu$ and $\pi$ using an independent sample, then constructing estimates of the unobserved outcomes and weights to be used in a second stage regression. This procedure is referred to as `quasi-oracle' because the error bound for the $\lambda$ estimator may decay faster than those for the $\mu$ and $\pi$ estimators, as has been shown by \cite{Nie2017} when penalized kernel ridge regression is used in the second stage. However, the R-learner framework remains `model-agnostic' in terms of the regression methods used in the first and second stage.

We propose a similar approach for learning $\beta$, which appears as an inverse weight in the estimator $\hat{\psi}$. Inverse weighting may be problematic when $\beta(x)$ is small, as small errors in $\hat{\beta}(x)$ can result in large differences in the value of $1/\hat{\beta}(x)$, a problem which is well documented in the context of inverse probability weighting estimators of the ATE \citep{Kang2007}. Concerns regarding extreme weights can be mitigated by regularizing the function $1/\hat{\beta}$ rather than $\hat{\beta}$, which can be achieved with reference to Lemma \ref{R_identLemma}. In particular, when $V=1$ and $W=\{A-\pi(X)\}^2$, the target ratio is $1/\beta(x)$. Therefore, an oracle estimator for $1/\beta$ is obtained by regressing $\{A-\pi(X)\}^{-2}$ on $X$ with weights $\{A-\pi(X)\}^2$. Analogously to the R-learner, we propose a learner that mimics this oracle learner by first estimating $\pi$ using an independent sample, then constructing estimates of the unobserved outcomes and weights to be used in a second stage regression.

\subsection{Proposed algorithms}
\label{proposed_algos}

The proposed working function estimators are implemented in Algorithms noSS and SS below. The latter uses a cross fitting regime to ensure that $\hat{\mu}(x_i),\hat{\pi}(x_i),\hat{\lambda}(x_i)$, and $\hat{\beta}(x_i)$ are obtained using working models which are constructed from a dataset that does not include the $i$th observation. This is useful in controlling the so-called empirical process term \citep{Chernozhukov2018,van_der_laan_cross-validated_2011}.
Algorithms noSS and SS return estimates $\{\hat{\pi}_i\}_{i=1}^n$, $\{\hat{\mu}_i\}_{i=1}^n$, $\{\hat{\lambda}_i\}_{i=1}^n$, and $\{\hat{\beta}_i\}_{i=1}^n$, which can be used to obtain
\begin{align*}
\hat{\psi} &= n^{-1}\sum_{i=1}^n  \frac{(a_i - \hat{\pi}_i)}{\hat{\beta}_i}\{y_i - \hat{\mu}_i - \hat{\lambda}_i(a_i - \hat{\pi}_i)\} +  \hat{\lambda}_i\\
\hat{\Psi} &= \frac{\sum_{i=1}^n(a_i - \hat{\pi}_i)(y_i - \hat{\mu}_i)}{\sum_{i=1}^n(a_i - \hat{\pi}_i)^2},
\end{align*}
with variances respectively estimated by $n^{-2}\sum_{i=1}^n \phi_{\psi,i}^2$ and $n^{-2}\sum_{i=1}^n \phi_{\Psi,i}^2$ where
\begin{align*}
\phi_{\psi,i} &= \frac{(a_i - \hat{\pi}_i)}{\hat{\beta}_i}\{y_i - \hat{\mu}_i - \hat{\lambda}_i(a_i - \hat{\pi}_i)\} +  \hat{\lambda}_i - \hat{\psi}\\
\phi_{\Psi,i} &= \frac{(a_i - \hat{\pi}_i)}{\hat{\eta}} \{y_i - \hat{\mu}_i - \hat{\Psi}(a_i - \hat{\pi}_i)\}
\end{align*}
and $\hat{\eta} \equiv n^{-1}\sum_{i=1}^n (a_i - \hat{\pi}_i)^2$ is an estimate of $\E\{\beta(X)\}$. We note that where the algorithms require regression estimates to be `fitted', any suitable regression / machine learning method can be used.

Both algorithms are also indexed by the choice of learner for $\hat{\lambda}$ and $\hat{\beta}$ in steps 2 and 3 of each algorithm respectively, with the substeps marked (A) and (B) referring to the direct, and quasi-oracle approaches. Note that the quasi-oracle methods in these algorithms do not use sample splitting to learn the unobserved outcomes and weights, this is due to the impracticality of the extra sample splitting required to do so. Substeps (A) and (B) do not need to be carried out for inference of $\hat{\Psi}$ only.

For estimators such as $\hat{\Psi}$, it has been suggested that faster convergence rates may be attained through additional sample splitting, which ensures that $\hat{\mu}(x_i)$ and $\hat{\pi}(x_i)$ are obtained from two different and independent datasets, neither of which containing the $i$th observation \citep{Newey2018}. We do not consider such `double cross fitting' here, since extensions, to estimate $\psi$, would require significant additional sample splitting to estimate $\{\hat{\lambda}_i\}_{i=1}^n$ and $\{\hat{\beta}_i\}_{i=1}^n$,  which may be impractical in finite samples.

\textbf{Algorithm noSS - no sample splitting}
\begin{enumerate}[label=(\arabic*)]
\item Fit $\hat{\mu}(x)$ and $\hat{\pi}(x)$. Use these fitted models to obtain $\hat{\mu}_i \equiv\hat{\mu}(x_i)$ and $\hat{\pi}_i \equiv\hat{\pi}(x_i)$.
\item (A) Fit $\Ehat\{YA|X=x\}$ and $\Ehat\{A^2|X=x\}$ and use these to construct $\hat{\lambda}(x)$ and $\hat{\beta}(x)$ as in \eqref{dir_lambda} and \eqref{dir_beta}. Or (B) obtain $\hat{\lambda}(x)$ and $1/\hat{\beta}(x)$ respectively by regressing $\{Y-\hat{\mu}(X)\}/\{A-\hat{\pi}(X)\}$ and $\{A-\hat{\pi}(X)\}^{-2}$ on $X$ with weights $\{A-\hat{\pi}(X)\}^2$ using all the data. After doing (A) or (B), use the fitted models to obtain $\hat{\lambda}_i \equiv\hat{\lambda}(x_i)$ and $\hat{\beta}_i \equiv\hat{\beta}(x_i)$.
\end{enumerate}

\textbf{Algorithm SS - sample splitting}
\begin{enumerate}[label=(\arabic*)]
\item Split the data into $K$ folds.
\item \textbf{For} each fold $k$: Fit $\hat{\mu}(x)$ and $\hat{\pi}(x)$ using the data set excluding fold $k$. Use these fitted models to obtain $\hat{\mu}_i \equiv\hat{\mu}(x_i)$ and $\hat{\pi}_i \equiv\hat{\pi}(x_i)$ for $i$ in fold $k$.
\item (A) Fit $\Ehat\{YA|X=x\}$ and $\Ehat\{A^2|X=x\}$ using the data set excluding fold $k$, and use these to construct $\hat{\lambda}(x)$ and $\hat{\beta}(x)$ as in \eqref{dir_lambda} and \eqref{dir_beta}. Or (B) obtain $\hat{\lambda}(x)$ and $1/\hat{\beta}(x)$ respectively by regressing $\{Y-\hat{\mu}(X)\}/\{A-\hat{\pi}(X)\}$ and $\{A-\hat{\pi}(X)\}^{-2}$ on $X$ with weights $\{A-\hat{\pi}(X)\}^2$ using the data set excluding fold $k$. After doing (A) or (B), use the fitted models to obtain $\hat{\lambda}_i \equiv\hat{\lambda}(x_i)$ and $\hat{\beta}_i \equiv\hat{\beta}(x_i)$ for $i$ in fold $k$. \textbf{End for.}
\end{enumerate}

\section{Simulation study}
\label{sect_simstud}

In our simulation study we compared Algorithms noSS and SS for estimating $\Psi$ and Algorithms noSS-A, noSS-B, SS-A and SS-B for estimating $\psi$ on generated data in finite samples, using $K=5$ fold sample splitting. Reproduction code for this study is available at \verb|github.com/ohines/alse|. We generated $1000$ datasets of size $n\in \{500,1000,...,4000\}$ from the following structural equation model
\begin{align*}
X_1,X_2,X_3 &\sim \text{Uniform}(-1,1) \\
\epsilon_1,\epsilon_2 &\sim \n{0}{1}\\
A &= X_1+0.5 X_1^3-2 X_2^2 + X_1^2 X_2 + (1+X_1^2)\epsilon_1\\
Y &= A(1+X_1-X_1^2-0.5 X_2^2) -X_1^2 X_2 + X_2 X_3  + \epsilon_2
\end{align*}
with the least squares estimands taking the true values $\psi = 0.5$ and $\Psi =107/294 \approx 0.36$. Note that, since $A$ is normally distributed given $X$, $\psi$ coincides with the unweighted ADE in this study.

\begin{figure}[htbp]
\centering
\includegraphics[width=\linewidth]{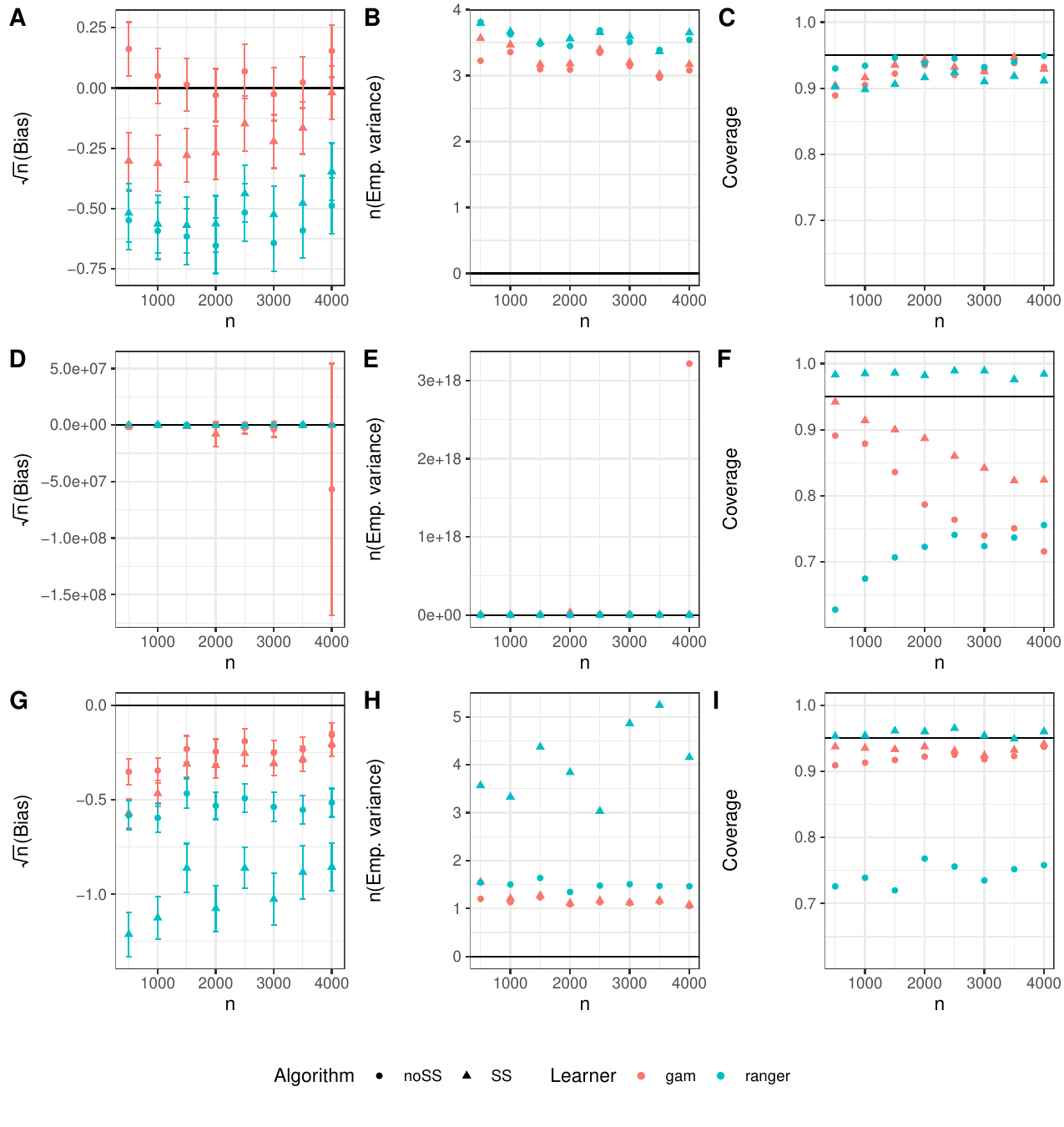}
\caption{Sample size against Bias (plots A,D,G), variance (plots B,E,H) and 95\% Wald CI coverage (plots C,F,I) for $\hat{\Psi}$ (plots A,B,C), $\hat{\psi}$ using the direct approach (plots D,E,F), and $\hat{\psi}$ using the quasi-oracle approach (plots G,H,I). We highlight that the y-axis limits change between rows of the bias and variance plots. The black horizontal lines represents zero bias, zero variance, and 95\% coverage respectively.
}
\label{sim_plot}
\end{figure}

For each dataset, $\hat{\psi}$ and $\hat{\Psi}$ were estimated along with their variance and associated Wald based (95\%) confidence intervals. Two regression model approaches were considered, the first used generalised additive models, as implemented through the \verb|mgcv| package in R \citep{Wood2016}. These models use flexible spline smoothing including pairwise interaction terms. The second regression modelling approach used random forest learners available through the \verb|ranger| package in R \citep{Wright2017}.
Figure \ref{sim_plot} shows empirical estimates of the bias and variance of $\hat{\psi}$ and $\hat{\Psi}$ scaled by $\sqrt{n}$ and $n$ respectively, as well as the empirical coverage probability of a Wald based 95\% confidence-interval. Comparing Algorithms noSS and SS for the estimation of $\hat{\psi}$, we notice that sample splitting generally improves confidence interval coverage.

Additionally, for estimation of $\psi$ the quasi-oracle approach (Algorithm -B) outperforms the direct approach (Algorithm -A) in terms of reduced bias, variance and improved CI coverage. This is achieved since the quasi-oracle learning controls the smoothness of $\hat{\lambda}$ and the inverse weights $1/\hat{\beta}$, whereas direct learning does not, leading to the possibility of extreme inverse weighting in the estimator. Indeed, the extreme bias values presented for Algorithm -A are due to extreme negative estimates of $1/\beta(X)$, which are corrected by the quasi-oracle estimator (Algorithm -B). On the basis of these results, we recommend Algorithm noSS-B for estimation of $\psi$ and Algorithm noSS for estimation of $\Psi$.

\section{Warfarin dose example}
\label{sect_appliedIll}

We illustrate the proposed estimators using the IWPC \citep{Cunningham2009} dataset, which has also been reanalysed several times in literature on dynamic treatment rule estimation \citep{Schulz2021,Wallace2018,Chen2016}. Reproduction code for this analysis is available at \verb|github.com/ohines/alse|. The data consists of $n=1732$ patients receiving Warfarin therapy, which is a commonly prescribed  anticoagulant used to treat thrombosis and thromboembolism. We consider least squares estimands for the effect of Warfarin dose $(A)$ on international normalised ratio (INR)  $(Y)$, which is a measure of blood clotting function, given 13 other patient characteristics $(X)$, including genetic data, as described in \cite{Cunningham2009}.

All models were fitted using the Super Learner of \citep{VanDerLaan2007}, an ensemble learning method implemented in the \verb|SuperLearner| package in R. This used 20 cross validation folds, and a `learner library' containing various routines (\verb|glm|, \verb|glmnet|, \verb|gam|, \verb|xgboost|, \verb|ranger|). Additional results using the `discrete' Super Learner for model fitting are presented in \myappendix \ref{appendix:iwpc_extra}. The discrete Super Learner selects the regression algorithm in the learner library which minimises a cross validated estimate of e.g. the mean squared error loss, whereas the Super Learner minimizes the same loss by taking a convex combination of learners. For the sample splitting algorithms (Algorithm SS), $K=20$ folds were chosen (between 10 to 20 folds is typical for cross-fitting procedures).

\begin{figure}[htb]
    \centering
    \includegraphics[width=0.9\linewidth]{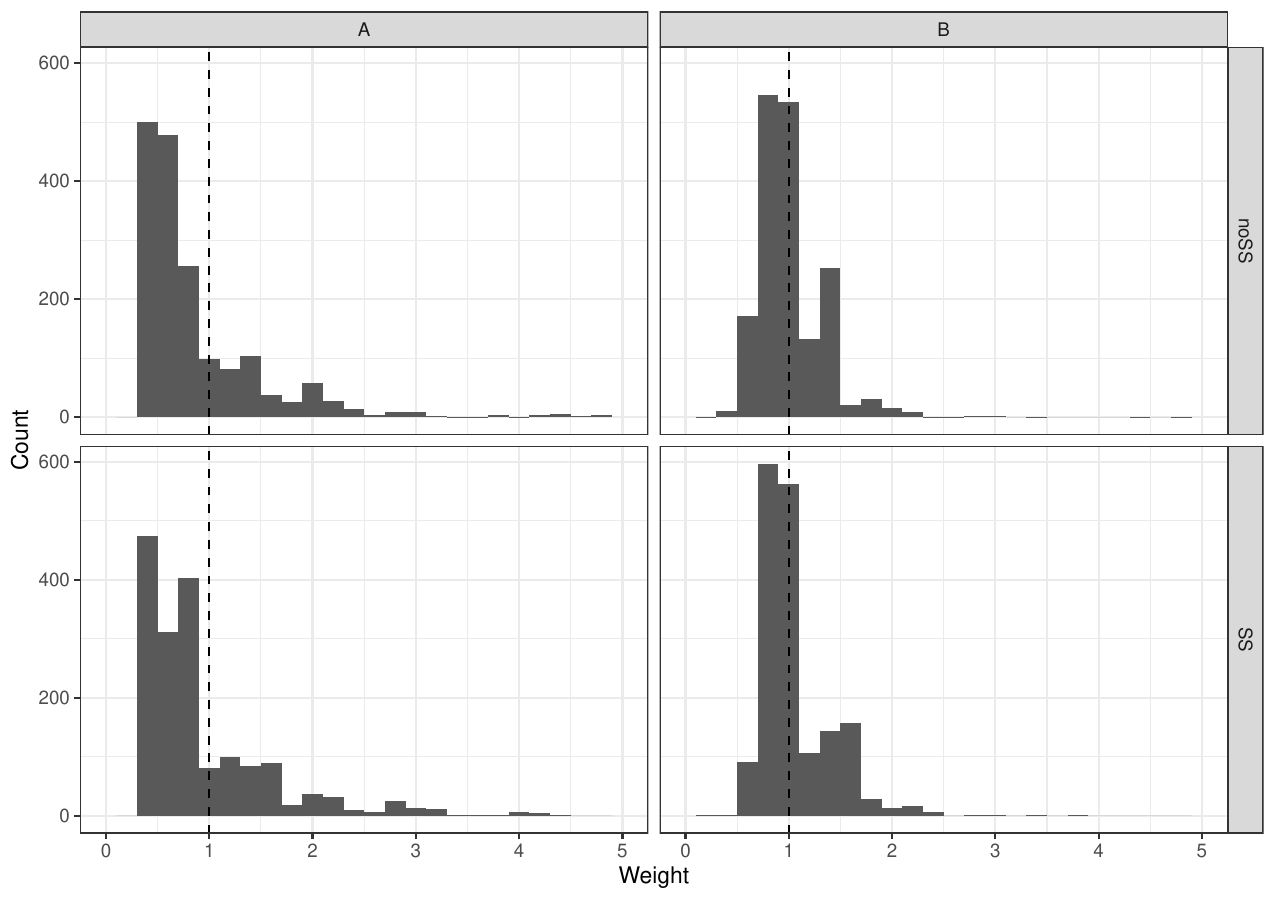}
    \caption{Least squares estimand weights approximated using the location-scale procedure described in \myappendix \ref{weight_approximations}. This procedure uses estimates of the conditional mean and variance of $A$ given $X$, which are obtained using the algorithms in Section \ref{proposed_algos} using the Super Learner for model fitting. Rows and columns refer to different algorithms as labeled.}
    \label{fig:noniscrete-superlearner-weights}
\end{figure}

First, we examine how the least squares estimand weights in \eqref{ALSE-weight-new} look in this example. These weights are one when $A$ is normal given $X$, in which case $\psi$ recovers the unweighted ADE. We approximate weights using the location-scale treatment model of \cite{Klyne2023} described in \myappendix \ref{weight_approximations}. Although this is not a necessary step for estimation of $\psi$ and $\Psi$, such approximations provide intuition as to how least squares estimands compare with unweighted average derivative estimands.
The histogram plots in Figure \ref{fig:noniscrete-superlearner-weights} suggest that $\psi$ and $\Psi$ can be interpreted as ADEs with a modest reweighting of the observations.

The results in Table \ref{IWPCresults} suggest that an increased Warfarin dose causes an increase in INR at a rate indicated by $\hat{\psi}$ and $\hat{\Psi}$. We see that the estimators for $\Psi$ have narrower confidence intervals than estimators for $\psi$, with commensurately smaller P-values for the corresponding `zero-effect' null ($\Psi = 0$ and $\psi=0$). This efficiency difference may be expected, since $\hat{\Psi}$ is based on the optimally efficient estimand in Section \ref{crump_type_result}, and since $\hat{\psi}$ requires additional nuisance function fitting, beyond that required for $\hat{\Psi}$. Additionally, the estimators for $\psi$, which use the R-learner (noSS-B and SS-B) give more credible estimates than those that use the direct approach (noSS-A and SS-A), in the sense that they are of a similar order of magnitude to the $\Psi$ estimates. Moreover, we see that sample splitting leads to more credible estimates, compared with no sample splitting, as evident in Algorithms SS-A versus noSS-A. This difference is likely because sample splitting helps to control for overfitting of the function estimators.

In this analysis, $\hat{\Psi}$ provides strong evidence for the causal effect of Warfarin dose on INR, but the better interpretability of $\psi$ means that $\hat{\psi}$ may provide better guidance as to the `average effect' of a small increase in Warfarin dose for the study population.

\begin{table}[htb]
\caption{Least squares estimands applied to IWPC data. Results indicate the point estimates, its standard error, and 95\% Wald confidence interval, all in units of INR/(mg/week). P-values correspond to a Wald test of the null hypothesis that the estimand is zero.}
\centering
\begin{tabular}{|l|l|l|l|l|l|}
\hline
Estimand            & Algorithm                 & Estimate & SE & CI & p\\ \hline
$\Psi$   & noSS      & 1.98$\times10^{-3}$ & 6.58$\times10^{-4}$ &(0.692$\times10^{-3}$, 3.27$\times10^{-3}$) & 0.003    \\ \hline
$\Psi$   & SS      & 1.89$\times10^{-3}$ & 6.26$\times10^{-4}$ & (0.662$\times10^{-3}$, 3.12$\times10^{-3}$) & 0.003   \\ \hline
$\psi$   & noSS-A     & 0.260 & 0.200 & (-0.132, 0.651) & 0.19     \\ \hline
$\psi$   & SS-A     & 0.0910 & 0.126 & (-0.155, 0.337) & 0.47   \\ \hline
$\psi$   & noSS-B     & 1.57$\times10^{-3}$ & 8.40$\times10^{-4}$ & (-0.0748$\times10^{-3}$, 3.22$\times10^{-3}$) & 0.06 \\ \hline
$\psi$   & SS-B     & 1.34$\times10^{-3}$ & 9.46$\times10^{-4}$ & (-0.510$\times10^{-3}$, 3.20$\times10^{-3}$) & 0.15   \\ \hline
\end{tabular}
\label{IWPCresults}
\end{table}

\section{Discussion}

We recommend that the least squares estimands $\psi$ and $\Psi$ be inferred for the exploratory analysis of the main effect of a continuous treatment on an outcome. The resulting quantities may be interpreted as weighted ADEs with data-adaptive weights. For $\psi$, these weights happen to equal $1$ when the treatment is conditionally normally distributed; and for $\Psi$ these weights happen to equal $1$ when the treatment is conditionally normally distributed and $\Var(A|X)$ is constant. For estimation, we recommend efficient one-step estimators, $\hat{\Psi}$ and $\hat{\psi}$, the latter of which generalises the AIPW to continuous treatments. To estimate working models we recommend a quasi-oracle approach based on the R-learner of \cite{Nie2017}, and a novel learner for the inverse variance, $1/\Var(A|X)$ that mitigates extreme weighting in the estimator.

%% file: supplement.tex
\section{Weighted average derivative effect proofs}
\label{appendix:optimal_derivs}

\subsection{Proof of weighted ADE identification}
\label{proof:identification}

Assume: (i) $Y^a\indep A \mid X$; (ii) $A=a \implies Y^a=Y$; (iii) $p_0(A|X)$ is continuous in $A$ almost surely; (iv) $\mu^\prime(a, x)$ exists.

Claim: Under these assumptions $\E\{w(A,X)\lim_{\varepsilon \to 0}\delta_{\varepsilon}(A,X)\} = \E\{w(A,X)\mu^\prime(A,X)\}$.

Proof: By (i) and (ii) we write
\begin{align*}
    \varepsilon \delta_{\varepsilon}(a,x) &=\E(Y^{a+\varepsilon}| X=x) - \E(Y^a | X=x) \\
    &= \E(Y^{a+\varepsilon}| A={a+\varepsilon}, X=x) - \E(Y^a | A=a, X=x) \\
    &= \E(Y | A={a+\varepsilon}, X=x) - \E(Y | A=a, X=x) \\
    &= \mu(a+\varepsilon, x) - \mu(a, x)
\end{align*}
where the extra conditioning in the second line is possible when $p_0(a|x) > 0$ and $p_0(a+\varepsilon|x) > 0$. If $p_0(a|x)$ is continuous in $a$, then $p_0(a+\varepsilon|x) > 0$ when $p_0(a|x) > 0$ and $\varepsilon$ is sufficiently small. Therefore, using (iii) and the fact that $p_0(A|X) > 0$ almost surely
\begin{align*}
    \lim_{\varepsilon \to 0} \delta_{\varepsilon}(A,X) &= \lim_{\varepsilon \to 0} \frac{\mu(A+\varepsilon, X) - \mu(A, X)}{\varepsilon} \\
    &= \mu^\prime(A, X)
\end{align*}
where each equality holds almost surely. Taking expectations we obtain
\begin{align*}
    \E\left\{w(A, X)  \lim_{\varepsilon \to 0} \delta_{\varepsilon}(A,X) \right\} &= \E\left\{w(A, X)  \mu^\prime(A, X) \right\}.
\end{align*}

\subsection{Proof of (\ref{canonical_contrast})}
\label{append:integration_by_parts}

Claim: Under (C1) and (C2), $\E\{w(A,X)f^\prime(A,X)|X\} = \E\{\alpha_w(A,X)f(A,X)|X\}$ almost surely for all differentiable $f\in \mathcal{H}$. Taking expectations, we further obtain $\E\{w(A,X)f^\prime(A,X)\} = \langle \alpha_w, f\rangle$, with $\alpha_w \in \mathcal{H}$ by (C3).

Proof: (C1) implies that $\Omega(x) = (s, t)$ where $s=s(x)$ and $t=t(x)$ are the boundaries of the support. Since the support is open, $p_0(s|x) = p_0(t|x) = 0$. Let $\tilde{p}(a,x)\equiv w(a,x)p_0(a|x)$ with $\tilde{p}(s,x) = \tilde{p}(t,x) = 0$. Using integration by parts,
\begin{align*}
\E\{w(A,X)f^\prime(A,X)|X = x\} &= \int_s^t  f^\prime(a,x) \tilde{p}(a,x) da \\
&=f(t,x)\tilde{p}(t,x)- f(s,x)\tilde{p}(s,x) - \int_s^t  f(a,x) \tilde{p}^\prime(a,x)  da \\
&= - \int_s^t  f(a,x) \left\{ w^\prime(a,x) p_0(a|x) + w(a,x) p_0^\prime(a|x)\right\} da \\
&= - \int_s^t  f(a,x) \left\{ w^\prime(a,x) + \frac{w(a,x) p_0^\prime(a|x)}{p_0(a|x)}\right\} p_0(a|x) da \\
&= \E\{ \alpha_w(A,X)f(A,X)|X = x\}
\end{align*}
where in the penultimate line we use the fact that $p_0(a|x) > 0$ for $a \in (s, t)$.

Remark: To motivate Theorem \ref{alpha_theorem}, note that the fundamental theorem of calculus implies
\begin{align*}
\tilde{p}(a,x) &= -\int_{s}^{a} \alpha_w(a^*,x)p_0(a^*| x) da^* .
\end{align*}

\subsection{Proof of Theorem \ref{alpha_theorem}}
\label{main_theorem_proof}

First note that the law of total expectation implies
\begin{align*}
    \E\{\alpha(A,X)|X\} &= \E\{\alpha(A,X)|A\leq a, X\}F(a|X) + \E\{\alpha(A,X)|A > a, X\}\{1 - F(a|X)\}.
\end{align*}
Since $\E\{\alpha(A,X)|X\} = 0$, for $\alpha \in \mathcal{R}$, we rewrite the weight in \eqref{f_tilde_theorem_alpha} as
\begin{align*}
    w(a,x)
    &= \frac{F(a|x)\{1-F(a|x)\}}{p_0(a|x)}\left[\E\{\alpha(A,X)|A > a, X=x\} - \E\{\alpha(A,X)|A \leq a, X=x\}\right] \\
    &= -\frac{\E\{\alpha(A,X)|A\leq a,X=x\}F(a|x) }{p_0(a|x)} \\
    &= \frac{-1}{p_0(a|x)}\int_{s}^{a} \alpha(a^*,x)p_0(a^*| x) da^* 
\end{align*}
where $s,t$ are the lower and upper boundary of the support of $A$. Therefore, for this weight: (C1) is satisfied by conditions of the Theorem; (C2) is satisfied by the fundamental theorem of calculus with \eqref{canonical_contrast} giving $\alpha_w = \alpha$;  and (C3) is satisfied because $\alpha \in \mathcal{H}$. 

For completeness, we explicitly show that $\E\{w(A,X)\} = 1$ for all $\alpha \in \mathcal{R}$, and $\E\{w(A,X)\mid X\} = 1$ almost surely for all $\alpha \in \mathcal{R}^*$. First write
\begin{align}
  w(a,x)p_0(a|x) &= - \int_s^t \alpha(a^*,x)\mathbb{I}(a^* \leq a) p_0(a^*|x) da^*  \label{weight_integral}
\end{align}
where $\mathbb{I}(a^* \leq a)$ is 1 when $a^* \leq a$ and 0 otherwise. It follows that
\begin{align*}
  \E\left\{w(A,X) | X = x \right\} &= \int_s^t w(a,x)p_0(a|x) da \\
  &= -\int_s^t \int_s^t \alpha(a^*,x)\mathbb{I}(a^* \leq a) p_0(a^*|x) da^* da \\
  &= \int_s^t \alpha(a^*,x) \left[\int_s^t -\mathbb{I}(a^* \leq a) da \right] p_0(a^*|x) da^* \\
  &= \int_s^t \alpha(a^*,x) \left[ a^* - t \right] p_0(a^*|x) da^* \\
  &= \E\{\alpha(A,X)A|X=x\} - t \E\{\alpha(A,X)|X = x\} \\
  &= \E\{\alpha(A,X)A|X=x\}.
\end{align*}
Hence $\E\left\{w(A,X) | X = x \right\} = 1$ for $\alpha \in \mathcal{R}^*$ and $\E\left\{w(A,X) \right\} = \E\left\{ \alpha(A,X)A \right\} = 1$ for $\alpha \in \mathcal{R}$.

\subsection{Proof of Theorem \ref{optimality_1}}
\label{optimality_proof}

Background: An efficient estimator $\hat{\theta}_w $ of $\theta_w $ is regular asymptotically linear, such that
\begin{align*}
\hat{\theta}_w = \theta_w + n^{-1}\sum_{i=1}^n \phi_{\theta,w}(o_i) + o_p(n^{-1/2})
\end{align*}
where $\phi_{\theta,w}(o)$ is the influence curve in \eqref{newey_ic}. Hence,
\begin{align*}
\sqrt{n}(\hat{\theta}_w - \theta_{w,S}) = n^{-1/2}\sum_{i=1}^n \alpha_w(a_i, x_i)\{y_i-\mu(a_i,x_i)\} + o_p(1)
\end{align*}
By the central limit theorem, $\sqrt{n}(\hat{\theta}_w - \theta_{w,S}) \overset{d}{\to}\n{0}{V}$, where the efficiency bound is
\begin{align*}
V &= \E\left\{\alpha_w^2(A,X)\{Y-\mu(A,X)\}^2\right\} \\
&= \E\left\{\alpha_w^2(A,X)\sigma^2(A,X)\right\} 
\end{align*}

Claim: $\E\left\{\alpha^2(A,X)\sigma^2(A,X)\right\}$ is minimised over $\alpha \in \mathcal{R}$ by the RR in Theorem \ref{optimality_1}.

Proof: Let $\gamma \in \mathcal{H}$ be a function with conditional mean $\bar{\gamma}(x) \equiv \E\{\gamma(A,X)|X = x\}$. The function $\alpha(a,x) = \gamma(a,x) - \bar{\gamma}(x)$ satisfies the condition $\E\{\alpha(A,X)|X\} = 0$, and hence, $\alpha \in \mathcal{R}$ if and only if $\E\{\alpha(A,X) A \} = \E[\alpha(A,X) \{A - \pi(X)\}] = 1$. We therefore consider the Lagrangian
\begin{align*}
  L(\gamma, \lambda) = \E\left( \{\gamma(A,X) - \bar{\gamma}(X)\}^2 \sigma^2(A, X) -2\lambda [\{\gamma(A,X) - \bar{\gamma}(X)\} \{A-\pi(X)\} -1 ] \right)
\end{align*}
where $\lambda$ is a Lagrange multiplier.
Next, consider replacing $\gamma(A,X)$ with an alternative $\gamma(A,X) + \delta \eta(A,X)$, where $\delta >0$ is constant. In doing so we obtain,
\begin{align*}
  \frac{L(\gamma + \delta \eta, \lambda) - L(\gamma, \lambda)}{\delta} = 2\E\left( \{\eta(A,X) - \bar{\eta}(X)\}\left[\{\gamma(A,X) - \bar{\gamma}(X)\}\sigma^2(A, X) -\lambda \{A-\pi(X)\}\right] \right) + O(\delta)
\end{align*}
where $\bar{\eta}(x) = \E\{\eta(A,X)|X=x\}$. We apply the identity
\begin{align*}
  \E\left(  \bar{\eta}(X)\left[\{\gamma(A,X) - \bar{\gamma}(X)\}\sigma^2(A, X) -\lambda \{A-\pi(X)\}\right] \right) &= \E\left[\eta(A,X) \Cov\{\gamma(A,X), \sigma^2(A,X)|X\}\right]
\end{align*}
To write
\begin{align*}
  \frac{L(\gamma + \delta \eta, \lambda) - L(\gamma, \lambda)}{\delta} = 2\E\left\{\eta(A,X) l(A,X, \gamma, \lambda) \right\} + O(\delta)
\end{align*}
where
\begin{align*}
  l(a,x, \gamma, \lambda) &= \{\gamma(a,x) - \bar{\gamma}(x)\}\sigma^2(a, x) -\lambda \{a-\pi(x)\} - \Cov\{\gamma(A,X), \sigma^2(A,X)|X=x\}.
\end{align*}
By the fundamental lemma of the calculus of variations, we require that $l(a,x, \gamma, \lambda) = 0$, and hence, writing $\alpha(a,x) = \gamma(a,x) - \bar{\gamma}(x)$ we require that
\begin{align*}
  \alpha(a,x)\sigma^2(a, x) - \E\left[\alpha(A,X)\sigma^2(A, X) | X=x \right] = \lambda \{a - \pi(x)\}
\end{align*}
This equality is satisfied by
\begin{align*}
  \alpha(a, x) &= \frac{\lambda \{a - \tilde{\pi}(x)\}}{\sigma^2(a, x)}.
\end{align*}
The constraint $\E\{\alpha(A,X) | X\} = 0$ sets
\begin{align*}
  \tilde{\pi}(x) \equiv \frac{\E\{A / \sigma^2(A,X) | X = x\}}{\E\{1 / \sigma^2(A,X) | X = x\}},
\end{align*}
and the constraint $\E\{\alpha(A,X) A\} = 1$ sets
\begin{align*}
  \lambda  &= \E\left[\frac{\{A - \tilde{\pi}(X)\}A}{\sigma^2(A, X)} \right]^{-1},
\end{align*}
which completes the proof.

\subsection{Proof to Theorem \ref{optimality_2}}
\label{optimality_proof2}

This proof is similar to that in \myappendix \ref{optimality_proof}. First, consider that $V = \E[\E\{\alpha(A|X)\sigma^2(A,X)|X\}]$. Hence, it is sufficient to minimize $\E\{\alpha(A|X)\sigma^2(A,X)|X=x\}$ for each $x$, with the constraints that $\E\{\alpha_w(A,X)A|X=x\}=1$ and $\E\{\alpha_w(A,X)|X=x\}=0$. For Lagrange multipliers $\lambda_1 = \lambda_1(x)$ and $\lambda_2 = \lambda_2(x)$, which are constant given $x$, we consider the Lagrangian
\begin{align*}
    L(\alpha, \lambda_1, \lambda_2, x) = \E\{l(A,X,\alpha, \lambda_1, \lambda_2)|X=x\}
\end{align*}
where
\begin{align*}
    l(a,x,\alpha, \lambda_1, \lambda_2) = \alpha(a,x)\sigma^2(a,x) -2\lambda_1(x)\alpha(a,x) - 2\lambda_2(x)\{\alpha(a,x)a - 1\}
\end{align*}
Differentiating $l(a,x,\alpha, \lambda_1, \lambda_2)$ with respect to $\alpha(a,x)$ and setting equal to zero gives
\begin{align*}
    \alpha(a,x) = \frac{\lambda_1(x) + \lambda_2(x)a}{\sigma^2(a,x)}.
\end{align*}
Since $\lambda_1$ and $\lambda_2$ are not yet defined we reparameterise $\alpha$ as
\begin{align*}
    \alpha(a,x) = \frac{\lambda(x)\{a - \tilde{\pi}(x)\}}{\sigma^2(a,x)}.
\end{align*}
where $\lambda(x)$ and $\tilde{\pi}(x)$ are as yet undefined.
The constraint $\E\{\alpha(A,X)|X=x\}=0$ sets $\tilde{\pi}(x)$ as in the main text, and the constraint $\E\{\alpha(A,X)A|X=x\}=1$ sets
\begin{align*}
  \lambda(x)  &= \E\left[\frac{\{A - \tilde{\pi}(X)\}A}{\sigma^2(A, X)} ~\Big|~ X=x\right]^{-1},
\end{align*}
which completes the proof.

\subsection{Theorem \ref{optimality_1} when the exposure is binary}
\label{crump_appendix}

Here we connect Theorem \ref{optimality_1} to Theorem 5.4 of \cite{Crump2006}, the latter of which considers the setting where $A \in \{0,1\}$ is a binary exposure. Letting $\sigma^2_a(x) = \sigma^2(a,x)$ we obtain
\begin{align*}
  \tilde{\pi}(x) &= \frac{\pi(x)\sigma^2_0(x)}{\pi(x)\sigma^2_0(x) + \{1 - \pi(x)\}\sigma^2_1(x)}.
\end{align*}
Hence, by checking the $a=1,0$ cases separately, one can verify that
\begin{align*}
  \frac{a - \tilde{\pi}(x)}{\sigma^2(a,x)} &= \frac{a - \pi(x)}{\pi(x)\sigma^2_0(x) + \{1 - \pi(x)\}\sigma^2_1(x)}.
\end{align*}
The unnormalised weight in \citet[Theorem 5.4]{Crump2006} is defined as
\begin{align*}
  w^*(x) &\equiv \left(\frac{\sigma^2_1(x)}{\pi(x)} + \frac{\sigma^2_0(x)}{1 - \pi(x)} \right)^{-1} \\
  &= \frac{\pi(x) \{1 - \pi(x)\}}{\pi(x)\sigma^2_0(x) + \{1 - \pi(x)\}\sigma^2_1(x)}
\end{align*}
Using this weight we write
\begin{align*}
  \frac{a - \tilde{\pi}(x)}{\sigma^2(a,x)} = \frac{w^*(x) \{a - \pi(x)\} }{\pi(x) \{1 - \pi(x)\}}
\end{align*}
Thus we obtain the normalisation constant
\begin{align*}
  \E\left\{ \frac{\{A - \tilde{\pi}(X)\} A }{\sigma^2(A,X) } \right\} &= \E\left\{ w^*(X) \frac{ \E[\{A - \pi(X)\} A| X] }{\pi(X) \{1 - \pi(X)\}} \right\} \\
  &= \E\{w^*(X)\}.
\end{align*}
Hence the RR in Theorem \ref{optimality_1} reduces to 
\begin{align*}
  \alpha_w(a,x) &= \left(\frac{w^*(x)}{\E\{w^*(X)\}} \right) \frac{\{a - \pi(x)\} }{\pi(x) \{1 - \pi(x)\}}
\end{align*}
which yields the optimally weighted ATE
\begin{align*}
  \langle \mu, \alpha_w \rangle &= \frac{\E[w^*(X) \{ \mu(1,X) - \mu(0,X)\}] }{\E\{w^*(X)\}}.
\end{align*}

\subsection{Proof of Lemma \ref{R_identLemma}}

Consider that
\begin{align*}
\E\left[ W\left\{\frac{V}{W} - g(X) \right\}^2 \right] &= \E\left\{ \frac{V^2}{W} \right\} + \E\left\{ g^2(X)W - 2g(X)V\right\} \\
&= \E\left\{ \frac{V^2}{W} \right\} + \E\left\{ g^2(X) \E(W|X) - 2g(X)\E(V|X)\right\}.
\end{align*}
Hence,
\begin{align*}
g^* &= \argmin_{g \in \mathcal{F}} \E\left[ W\left\{\frac{V}{W} - g(X) \right\}^2 \right] \\
&= \argmin_{g \in \mathcal{F}} \E\left\{ g^2(X) \E(W|X) - 2g(X)\E(V|X)\right\}.
\end{align*}
By the calculus of variations, $g^*(x)\E(W|X=x)  = \E(V|X=x)$. The result follows since $W>0$ almost surely implies that $\E(W|X) \neq 0$ almost surely.

\section{Investigation into the Least squares estimand weight}

\subsection{Least squares projection}
\label{sect:projection}

Here we describe how least squares estimands (Example \ref{alse}) are connected to model projections. For a model $\mathcal{M} \subseteq \mathcal{H}$ we define the least squares projection as the function $\tilde{\mu}_{\mathcal{M}} \in \mathcal{M}$ that is closest to the unknown regression function $\mu \in \mathcal{H}$ in the sense that
\begin{align*}
\tilde{\mu}_{\mathcal{M}} &\equiv \argmin_{f\in \mathcal{M}} \E \left[\{\mu(A,X)-f(A,X)\}^2 \right].
\end{align*}
This notion of model projection is considered by \cite{Neugebauer2007} and \cite{Chambaz2012}, who propose projections on to linear working models, and by \cite{Buja2019} who consider likelihood-based projections. Following \cite{Robinson1988}, we consider a semi-parametric partially linear model $\mathcal{M}_i$, which is the set of functions of the form $\omega(x) + \nu a$, indexed by the infinite dimensional parameter $(\omega, \nu)$, where $\omega: \mathbb{R}^p \to \mathbb{R}$ is a function and $\nu \in \mathbb{R}$ is a scalar. The projection of $\mu$ on to $\mathcal{M}_i$ is $\tilde{\mu}_{\mathcal{M}_i}(a,x) = \mu(x) + \Psi \{a-\pi(x)\}$.
Hence we say that $\Psi$ is a `least squares estimand' as it is the coefficient in a partially linear projection model which minimises the mean squared remainder. Crucially, the model $\mathcal{M}_i$ is used to interpret the nonparametrically defined estimand $\Psi$, but we do not assume that the model is `true', in the sense that we do not require that $\mu \in \mathcal{M}_i$. 

This projection view of least squares estimands is extended by considering the more flexible conditionally linear model $\mathcal{M}_{j} \supset \mathcal{M}_i$, which is the set of functions of the form $\omega(x) + \nu(x) a$, indexed by the functions $\omega$ and $\nu: \mathbb{R}^p \to \mathbb{R}$. The projection of $\mu$ on $\mathcal{M}_j$ is $\tilde{\mu}_{\mathcal{M}_j}(a,x) = \mu(x) + \lambda(x) \{a-\pi(x)\}$.
Hence $\psi$ and $\Psi$ are `least squares estimands' in the sense that they are weighted averages of the conditional least squares function $\lambda(x)$, i.e. they are of the form $\E\{w(X)\lambda(X)\}$, where $w(x) = 1$ for $\psi$ and $w(x) = \beta(x) / \E\{\beta(X)\}$ for $\Psi$, with $\E\{w(X)\} = 1$ in both cases. Note that, when $A$ is a binary treatment, $\mu \in\mathcal{M}_j$ and $\lambda(x)=\mu(1,x)-\mu(0,x)$ identifies the conditional ATE $\E(Y^1 - Y^0\mid X=x)$ when $X$ is sufficient to adjust for confounding.

The estimand $\Psi$ appears in \cite{Vansteelandt2020} who consider inference for the constant term indexing $h\{\mu(a,x)\}\in\mathcal{M}_i$ where $h: \mathbb{R} \to \mathbb{R}$ is a canonical link function. Rather than appealing to model projection, they set out desiderata for an estimand under model misspecification, with their resulting estimand reducing to $\Psi$ when $h$ is an identity link. Similarly, $\Psi$ appears in other work on partially linear models without reference to projection \citep{Newey2018,Robins2008}.
The fact that least squares estimands are weighted ADEs is a novel contribution of this work that is closely related to two observations in the literature:
\cite{Banerjee2007} construct an estimator of the ADE by partitioning the support of $X$ into disjoint bins, and applying a linear regression to each bin---the ADE estimate is a weighted average of regression coefficients in each bin;
% \cite{Buja2019} interpret ordinary least squares coefficients as weighted sums of `slopes' between pairs of observations, without invoking differentiability;
and \cite{Hirshberg2018} show that $\psi$ recovers the ADE when $\mu \in\mathcal{M}_j$. The key difference between the latter work and ours is that we do not rely on a functional form for $\mu$, rather we show that $\psi$ is an ADE with a certain kind of weighting.

\subsection{Parametric distributions}
\label{append_ls_weight}

Here we examine the ADE weight $w_\psi(a, x)$ in Example \ref{alse}, which is associated with the least squares estimands, $\psi$ and proportional to the weight $w_\Psi(a,x)$ for $\Psi$ in the main text. We remark that this weight depends on the unknown distribution $P_0$ through the distribution of $A$ given $X$. It is, however, not necessary to compute this weight to estimate $\psi$ (or $\Psi$).

We consider the form of this exposure weight under various parametric distributions $p_0(a|x)$ in Pearson's distribution family. Table \ref{weight_table} summarises the functional form of the weight, with details provided in Examples \ref{norm_ex} to \ref{t_ex} below. For each distribution in Table \ref{weight_table}, we report the ``conditionally unnormalised'' weight, which recovers the true weight as $w_\psi(a, x) = \tilde{w}(a, x) / \E\{\tilde{w}(A, X)|X=x\}$.

\begin{table}[htbp]
\caption{Conditionally unnormalised weight $\tilde{w}(a, x)$ for least squares estimands when the distribution of $A$ given $X$ follows a parametric distribution. Proofs below.}
\label{weight_table}
\centering
\begin{tabular}{|l|l|l|}
\hline
Exposure Distribution & Support & $\tilde{w}(a, x)$ \\ \hline
Normal & $A\in(-\infty, \infty)$ & $1$ \\
Gamma & $A>0$ & $a$ \\
Inverse-Gamma & $A>0$ & $a^2$ \\
Beta & $A\in(0,1)$ & $a(1-a)$ \\
Beta Prime & $A>0$ & $a(1+a)$ \\
t-distribution (d.o.f. $=\nu(x)$) & $A\in(-\infty, \infty)$ & $1+\frac{a^2}{\nu(x)}$ \\ \hline
\end{tabular}
\end{table}

When deriving the results in Table \ref{weight_table}, it is helpful to consider the function $\tilde{p}(a|x) \equiv w_\psi(a, x)p_0(a|x)$. Since the weight is non-negative (by Lemma \ref{sufficiency}), and $\E\{w_\psi(A, X)|X\} = 1$, it follows that $\tilde{p}(a|x)$ is a density function. Letting $s$ denote the lower boundary of the support of $A$ then we write this density function as
\begin{align*}
    \tilde{p}(a|x) = \int_{s}^{a} \frac{\mu(x) - a^*}{\beta(x)} p_0(a^*|x) da^*
\end{align*}
We use this function as a tool to derive the exposure weight in the following examples. Each example requires verifying a derivative result, which follows from standard calculus and the Gamma function property $\nu\Gamma(\nu) = \Gamma(\nu + 1)$.

\begin{example}[Normal distribution]
\label{norm_ex}
Letting $p_0(a|\mu, \beta)$ denote the normal distribution, with mean $\mu = \mu(x)$ and variance $\beta = \beta(x)$, we claim that $\tilde{p}(a|\mu, \beta) = p_0(a|\mu, \beta)$, which implies a weight $w(a, x) = \tilde{p}(a|\mu, \beta) / p_0(a|\mu, \beta) = 1$. To verify this claim, note that
\begin{align*}
\frac{d}{da} p_0(a|\mu, \beta) = \frac{\mu-a}{\beta} p_0(a|\mu, \beta).
\end{align*}
Hence, 
\begin{align*}
    \tilde{p}(a|x) = \int_{-\infty}^{a} \left\{\frac{d}{da^*} p_0(a^*|\mu, \beta)\right\} da^*
\end{align*}
and the result follows by the fundamental theorem of calculus. The unitary weight of the normal distribution implies that the least squares estimand $\psi$ recovers the average derivative effect when the exposure $A$ is normally distributed given $X$.
\end{example}

\begin{example}[Gamma distribution]
\label{gamma_ex}
The Gamma distribution, with shape parameter, $\nu = \nu(x) > 0$, and rate parameter $\tau = \tau(x) > 0$, has the density
\begin{align*}
p_0(a|\nu,\tau) = \frac{\tau^\nu}{\Gamma(\nu)} a^{\nu-1}\exp(-\tau a)
\end{align*}
for $a > 0$ and 0 otherwise. 

We claim that, for the gamma distribution, $\tilde{p}(a|\nu,\tau) = p_0(a|\nu+1,\tau)$. As in Example \ref{norm_ex}, it is sufficient to verify that
\begin{align*}
\frac{d}{da} p_0(a|\nu+1,\tau) = \frac{\mu-a}{\beta} p_0(a|\nu,\tau)
\end{align*}
where the mean and variance are $\mu=\nu/\tau$ and $\beta = \nu/\tau^2$. Therefore, the weight is
\begin{align*}
    w_\psi(a, x) = \frac{p_0(a|\nu+1,\tau)}{p_0(a|\nu,\tau)} = \frac{a}{\mu}.
\end{align*}
This linear weight assigns unitary weight to the mean value $a=\mu$, with larger weight given to values above the mean.
\end{example}

\begin{example}[Inverse Gamma distribution]
\label{inv_gamma_ex}
The Inverse Gamma distribution, with shape parameter, $\nu = \nu(x) > 0$, and scale parameter $\tau = \tau(x) > 0$, has the density
\begin{align*}
p_0(a|\nu,\tau) = \frac{\tau^\nu}{\Gamma(\nu)} a^{-\nu-1}\exp(-\tau / a)
\end{align*}
for $a > 0$ and 0 otherwise. We claim that, for the inverse-gamma distribution, $\tilde{p}(a|\nu,\tau) = p_0(a|\nu-2,\tau)$. As before, it is sufficient to verify that
\begin{align*}
\frac{d}{da} p_0(a|\nu-2,\tau) = \frac{\mu-a}{\beta} p_0(a|\nu,\tau)
\end{align*}
where the mean and variance are $\mu = \tau/(\nu-1)$ and
\begin{align*}
    \beta = \frac{\tau^2}{(\nu-1)^{2}(\nu-2)}.
\end{align*}
Therefore, the weight is
\begin{align*}
    w_\psi(a, x) = \frac{p_0(a|\nu-2,\tau)}{p_0(a|\nu,\tau)} = \frac{a^2}{\beta+\mu^2}.
\end{align*}
\end{example}

\begin{example}[Beta distribution]
\label{beta_ex}
The beta distribution, with shape parameters, $\nu = \nu(x) > 0$, and $\tau = \tau(x) > 0$, has the density
\begin{align*}
p_0(a|\nu,\tau) = \frac{\Gamma(\nu + \tau)}{\Gamma(\nu)\Gamma(\tau)} a^{\nu-1}(1-a)^{\tau-1}
\end{align*}
for $a \in [0,1]$ and 0 otherwise. We claim that, for the beta distribution, $\tilde{p}(a|\nu,\tau) = p_0(a|\nu+1,\tau+1)$. As before, it is sufficient to verify that
\begin{align*}
\frac{d}{da} p_0(a|\nu+1,\tau+1) = \frac{\mu-a}{\beta} p_0(a|\nu,\tau)
\end{align*}
where the mean and variance are $\mu = \nu / (\nu + \tau)$ and
\begin{align*}
    \beta = \frac{\nu\tau}{(\nu+\tau)^{2}(\nu+\tau+1)}.
\end{align*}
Therefore, the weight is
\begin{align*}
    w_\psi(a, x) = \frac{p_0(a|\nu+1,\tau+1)}{p_0(a|\nu,\tau)} = \frac{a(1-a)}{\mu(1-\mu) - \beta}.
\end{align*}
This quadratic weight is at a maximum when $a=1/2$, with very little weight assigned to values close to $a=0,1$.
\end{example}

\begin{example}[Beta-prime distribution]
\label{beta_prime_ex}
The beta-prime distribution, with shape parameters, $\nu = \nu(x) > 2$, and $\tau = \tau(x) > 0$, has the density
\begin{align*}
p_0(a|\nu,\tau) = \frac{\Gamma(\nu + \tau)}{\Gamma(\nu)\Gamma(\tau)} a^{\nu-1}(1+a)^{-\nu-\tau}
\end{align*}
for $a \geq 0 $. We claim that, for the beta-prime distribution, $\tilde{p}(a|\nu,\tau) = p_0(a|\nu+1,\tau-2)$. As before, it is sufficient to verify that
\begin{align*}
\frac{d}{da} p_0(a|\nu+1,\tau-2) = \frac{\mu-a}{\beta} p_0(a|\nu,\tau)
\end{align*}
where the mean and variance are $\mu = \nu / (\tau - 1)$ and 
\begin{align*}
    \beta = \frac{\nu(\nu + \tau - 1)}{(\tau - 2)(\tau - 1)^2}.
\end{align*}
Therefore, the weight is
\begin{align*}
    w_\psi(a, x) = \frac{p_0(a|\nu+1,\tau-2)}{p_0(a|\nu,\tau)} = \frac{a(1+a)}{\mu(1+\mu) + \beta}.
\end{align*}
%Like the Gamma exposure weight in Example \ref{gamma_ex}, the Beta-prime exposure weight is strictly increasing, however the Beta-prime weight increases quadratically, whereas the Gamma exposure weight increases linearly.

\end{example}

\begin{example}[Student's t-distribution]
\label{t_ex}
The t-distribution, with $\nu = \nu(x) > 2$ degrees of freedom, location parameter $\mu=\mu(x)$, and scale parameter $\tau = \tau(x)$ has the density
\begin{align*}
p_0(a|\nu, \mu, \tau^2) = \frac{\Gamma\left(\frac{\nu + 1}{2}\right)}{\sqrt{\nu\pi}\tau\Gamma\left(\frac{\nu}{2}\right)} \left(1 + \frac{(a-\mu)^2}{\nu \tau^2}\right)^{-\frac{\nu + 1}{2}}
\end{align*}
We claim that, for the t-distribution,
\begin{align*}
\tilde{p}(a|\nu, \mu, \tau^2) = p_0(a|\nu - 2, \mu, \beta)
\end{align*}
where $\beta = \tau^2\nu/(\nu-2)$ is the variance of the t-distribution. As before, it is sufficient to verify that
\begin{align*}
\frac{d}{da} p_0(a|\nu - 2, \mu, \beta) = \frac{\mu-a}{\beta} p_0(a|\nu, \mu, \tau^2).
\end{align*}
Therefore, the weight is
\begin{align*}
    w_\psi(a, x) = \frac{p_0(a|\nu - 2, \mu, \beta)}{p_0(a|\nu, \mu, \tau^2)} = \frac{\left(1+\frac{(a-\mu)^2}{\nu\tau^2}\right)}{\left(1+\frac{1}{\nu-2}\right)}.
\end{align*}
As $\nu \to \infty$ then $w(a|x) \to 1$. This is expected since the the t-distribution tends to a normal distribution, with mean $\mu$, and variance $\tau^2$ in this limit. Note that in Table \ref{weight_table}, this weight is reported for $\mu=0, \tau=1$.
\end{example}

\subsection{Numerical approximations}
\label{weight_approximations}

Here we propose a simple method to estimate least squares weights numerically. We emphasise that estimating weights is not necessary for inference or interpretation of least squares estimands, but may be useful for visually inspecting how, for a given dataset, least squares estimands might compare with an average derivative estimands, i.e. where $w(a,x) = 1$. Following \cite{Kennedy2017} and \cite{Klyne2023}, we assume a location-scale exposure model $A = \pi(X) + \beta^{1/2}(X) U$ where $U$ is a random variable with $U \indep X$. Under this model, we write the density
\begin{align*}
    p_0(a|x) &= p_u \left\{ u(a,x)\right\} \\
    u(a,x) &\equiv \frac{a - \pi(x)}{\beta^{1/2}(x)}
\end{align*}
where $p_u$ is the marginal density of $U$. Applying this density expression to the integral express of the least squares estimand weight in \eqref{weight_integral}, we obtain the weight
\begin{equation}
    w(a, x) = \frac{-1}{p_u \{ u(a,x) \}} \int_{-\infty}^\infty u^* \mathbb{I}\{u^* \leq u(a,x) \} p_u(u^*) du^* \label{location-scale}
\end{equation}
where we have replaced the bounds of the support of $A$ with $\pm \infty$, and used the fact that, for least squares estimands, $\alpha(a,x) =  u(a,x) \beta^{-1/2}(x)$. To approximate \eqref{location-scale} numerically, we first estimate $u(a,x)$ by centering and scaling the observed exposure using estimates $\hat{\pi}$ and $\hat{\beta}$ that are obtained from the algorithms described in Section \ref{proposed_algos}. Next we estimate $p_u$ using a kernel density estimator, and approximate the integral in \eqref{location-scale} empirically using Monte-Carlo.

\section{Estimator Asymptotic Distribution}
\label{appendix:asymptotics}

In this Appendix we use a common empirical processes notation, where we define linear operators $P$ and $\mathbb{P}_n$ such that for some function $h(O)$, $P\{ h(O)\} \equiv \E \{h(O)\}$ and $\mathbb{P}_n \{h(O)\} \equiv n^{-1} \sum_{i=1}^n h(o_i)$.

\subsection{Proof of Theorem \ref{asym_theorem_alse1}}
Define
\begin{align*}
\hat{\phi}_{\psi}(o) &= \frac{\{y-\hat{\mu}(x)\}\{a-\hat{\pi}(x)\}-\hat{\lambda}(x)\{a-\hat{\pi}(x)\}^2}{\hat{\beta}(x)} + \hat{\lambda}(x) -\hat{\psi}_0
\end{align*}
where $\hat{\psi}_0 = \mathbb{P}_n\{\hat{\lambda}(X)\}$ denotes an initial estimate of $\psi$. Without making any restrictions we write
\begin{align}
    \hat{\psi} - \psi &= (\mathbb{P}_n -P)\{\phi_{\psi}(O)\} + R_n + H_n \\
    \hat{\psi} &\equiv \hat{\psi}_0 + \mathbb{P}_n\{\hat{\phi}_{\psi}(O)\} \\
    R_n &\equiv \hat{\psi}_0 - \psi + P\{\hat{\phi}_{\psi}(O)\} \\
    H_n &\equiv (\mathbb{P}_n -P)\{\hat{\phi}_{\psi}(O) - \phi_{\psi}(O)\}.
\end{align}
We will show that the remainder therm $R_n=o_P(n^{-1/2})$ and the empirical process term $H_n=o_P(n^{-1/2})$, and hence the result follows since $P\{\phi_{\psi}(O)\} = 0$. 

\textbf{The remainder term}

To simplify notation, we will omit function arguments, e.g. $\mu = \mu(X)$ with similar for $\hat{\mu},\lambda,\hat{\lambda},\pi,\hat{\pi},\hat{\beta},\beta$. Since $\psi = \E[\lambda]$ we can write the remainder term as
\begin{align*}
R_n &= \E\left[ \frac{(Y-\hat{\mu})(A-\hat{\pi})-\hat{\lambda}(A-\hat{\pi})^2}{\hat{\beta}} + \hat{\lambda}-\lambda \right]
\end{align*}
Note that
\begin{align*}
\E[(Y-\hat{\mu})(A-\hat{\pi})|X] &= \lambda\beta + (\mu-\hat{\mu})(\pi-\hat{\pi}) \\
\E[(A-\hat{\pi})^2|X] &= \beta + (\pi-\hat{\pi})^2
\end{align*}
and hence
\begin{align*}
R_n &= \E\left[ \frac{(\pi-\hat{\pi}) \{\mu-\hat{\mu}-\hat{\lambda}(\pi-\hat{\pi}) \} + (\lambda-\hat{\lambda})(\beta-\hat{\beta}) }{\hat{\beta}} \right] 
\end{align*}
Using the inequality, $(a+b)^2\leq 2(a^2+b^2)$,
\begin{align}
\frac{R_n^2}{2} &\leq \underbrace{\E\left[ \frac{(\pi-\hat{\pi}) \{\mu-\hat{\mu}-\hat{\lambda}(\pi-\hat{\pi}) \}}{\hat{\beta}} \right]^2}_{\text{first remainder}}  + \underbrace{\E\left[\frac{(\lambda-\hat{\lambda})(\beta-\hat{\beta}) }{\hat{\beta}} \right]^2 }_{\text{second remainder}} \label{appendix_remainder}
\end{align}
We will show that the two remainder terms on the right hand side are $o_p(n^{-1})$. For the first remainder, the Cauchy-Schwarz inequality gives
\begin{align*}
\E\left[ \frac{(\pi-\hat{\pi}) \{\mu-\hat{\mu}-\hat{\lambda}(\pi-\hat{\pi}) \}}{\hat{\beta}} \right]^2 &\leq \E\left[\frac{(\pi-\hat{\pi})^2}{\hat{\beta}}\right] \E\left[\frac{\{\mu-\hat{\mu}-\hat{\lambda}(\pi-\hat{\pi}) \}^2}{\hat{\beta}}\right] \\
&\leq \left(\frac{1}{\epsilon^2}\right) \E\left[(\pi-\hat{\pi})^2\right] \E\left[\{\mu-\hat{\mu}-\hat{\lambda}(\pi-\hat{\pi}) \}^2\right] \\
&\leq \left(\frac{2}{\epsilon^2}\right) \E\left[(\pi-\hat{\pi})^2\right] \left\{\E\left[(\mu-\hat{\mu})^2\right] + \E\left[\hat{\lambda}^2(\pi-\hat{\pi})^2\right] \right\} \\
&\leq \left(\frac{2}{\epsilon^2}\right) \E\left[(\pi-\hat{\pi})^2\right] \left\{\E\left[(\mu-\hat{\mu})^2\right] + K^2\E\left[(\pi-\hat{\pi})^2\right] \right\}
\end{align*}
To obtain the second inequality above, we choose $\epsilon > 0$ such that $\hat{\beta} \geq \epsilon$ almost surely, and to obtain the third inequality we once again apply the inequality $(a+b)^2\leq 2(a^2+b^2)$. It therefore follows from (A1) that the first remainder in \eqref{appendix_remainder} is $o_p(n^{-1})$.

For the second remainder in \eqref{appendix_remainder}, the Cauchy-Schwarz inequality gives
\begin{align*}
\E\left[\frac{(\lambda-\hat{\lambda})(\beta-\hat{\beta})}{\hat{\beta}}  \right]^2 &\leq \E\left[\frac{(\lambda-\hat{\lambda})^2}{\hat{\beta}}\right] \E\left[\frac{(\beta-\hat{\beta})^2}{\hat{\beta}}\right] \\
&\leq \left(\frac{1}{\epsilon^2}\right) \E\left[(\lambda-\hat{\lambda})^2\right] \E\left[(\beta-\hat{\beta})^2\right]
\end{align*}
which is $o_P(n^{-1})$ under (A2). Hence $R_n=o_P(n^{-1/2})$.

\textbf{The empirical process term}

First write the empirical process term as the sum of four terms
\begin{align*}
    H_n &= (\mathbb{P}_n - P)\{\psi - \hat{\psi}_0\} \\
    &+ (\mathbb{P}_n - P)\{\hat{\lambda} - \lambda\} \\
    &+ (\mathbb{P}_n - P)\left\{\frac{\hat{u}(O)}{\hat{\beta}} - \frac{u(O)}{\beta}\right\} \\
    &+ (\mathbb{P}_n - P)\left\{\frac{\hat{\lambda}\hat{v}(O)}{\hat{\beta}} - \frac{\lambda v(O)}{\beta}\right\}
\end{align*}
where
\begin{align*}
    \hat{u}(O) &\equiv (Y-\hat{\mu})(A-\hat{\pi}) \\
    u(O) &\equiv (Y-\mu)(A-\pi) \\
    \hat{v}(O) &\equiv (A-\hat{\pi})^2 \\
    v(O) &\equiv (A-\pi)^2
\end{align*}
Note that the first term is zero since $(\mathbb{P}_n - P)\{\psi - \hat{\psi}_0\} = (\psi - \hat{\psi}_0)(\mathbb{P}_n - P)\{1\} = 0$. When the Donsker condition holds, then, by Lemma 19.24 of \cite{Vaart2013}, the second term is $o_P(n^{-1/2})$ provided (i) that $\E\{(\hat{\lambda} - \lambda)^2\} = o_p(1)$, the third term is $o_P(n^{-1/2})$ provided (ii) that $\E\left[\left\{\frac{\hat{u}(O)}{\hat{\beta}} - \frac{u(O)}{\beta}\right\}^2\right] = o_p(1)$, and the fourth term is $o_P(n^{-1/2})$ provided (iii) that $\E\left[\left\{\frac{\hat{\lambda}\hat{v}(O)}{\hat{\beta}} - \frac{\lambda v(O)}{\beta}\right\}^2\right] = o_p(1)$. Similarly, under sample splitting then by Chebyshev's inequality, (i), (ii), and (iii) are also sufficient conditions for $H_n$ to be $o_P(n^{-1/2})$. Note that condition (i) holds since $\hat{\lambda}(x)$ is a consistent estimator of $\lambda(x)$, therefore, we must only show that (ii) and (iii) hold.

For (ii) we write
\begin{align*}
    \frac{\hat{u}(O)}{\hat{\beta}} - \frac{u(O)}{\beta} &= \frac{\hat{u}(O) - u(O)}{\hat{\beta}} + \left\{\frac{1}{\hat{\beta}} - \frac{1}{\beta}\right\}u(O)  \\
    &= \frac{(Y-\mu) \epsilon_\pi + (A-\pi) \epsilon_\mu + \epsilon_\mu \epsilon_\pi}{\hat{\beta}} - \frac{\epsilon_\beta (Y-\mu)(A-\pi)}{\beta\hat{\beta}}
\end{align*}
where $\epsilon_\mu = \hat{\mu} - \mu$, $\epsilon_\pi = \hat{\pi} - \pi$, and $\epsilon_{\beta} = \hat{\beta} - \beta$. Hence, conditioning on the sample that delivered the nuisance parameter estimators (which we do not make explicit in our notation)
\begin{align*}
    \E\left[\left\{\frac{\hat{u}(O)}{\hat{\beta}} - \frac{u(O)}{\beta} \right\}^2 \Big| X\right] &= \frac{k_{2,0}(X)\epsilon_\pi^2 + k_{0,2}(X)\epsilon_\mu^2 + 2k_{1,1}(X)\epsilon_\mu\epsilon_\pi + \epsilon^2_\mu \epsilon_\pi^2}{\hat{\beta}^2} \\
    &+ \frac{\epsilon_{\beta}^2k_{2,2}(X)}{\beta^2\hat{\beta}^2} \\
    &+ \frac{2\epsilon_\beta}{\beta\hat{\beta}^2}\left\{k_{2,1}(X)\epsilon_\pi + k_{1,2}(X) \epsilon_\mu + k_{1,1}(X)\epsilon_\pi \epsilon_\mu \right\}
\end{align*}
where $k_{i,j}(X) = \E\{(Y-\mu)^i(A-\pi)^j|X\}$. By Cauchy-Schwarz, $k_{i,j}^2(X) \leq E\{(Y-\mu)^{2i}|X\}E\{(A-\pi)^{2j}|X\}$ hence each of the $k_{i,j}^2(X) \leq K^2$ in the expression above
\begin{align*}
   \Bigg| \E\left[\left\{\frac{\hat{u}(O)}{\hat{\beta}} - \frac{u(O)}{\beta} \right\}^2 \Big| X\right] \Bigg| \leq \frac{K}{\epsilon^2} \{\epsilon_{\pi}^2 + \epsilon_{\mu}^2 + 2 |\epsilon_\mu \epsilon_\pi|\} + \frac{\epsilon_{\pi}^2 \epsilon_{\mu}^2}{\epsilon^2} + \frac{K}{\epsilon^4} \epsilon_{\beta}^2 + \frac{2K}{\epsilon^3}|\{\epsilon_{\pi} + \epsilon_{\mu} + \epsilon_{\pi}\epsilon_{\mu}\}|
\end{align*}
Since each of $\pi,\mu,\beta$ are consistent, the right hand side above tends to zero, as does (ii) by the dominated convergence theorem. The proof of (iii) proceeds in a similar way,
\begin{align*}
    \frac{\hat{\lambda}\hat{v}(O)}{\hat{\beta}} - \frac{\lambda v(O)}{\beta} &= \frac{\hat{\lambda} \epsilon_{\pi} \{2(A-\pi) + \epsilon_{\pi} \} }{\hat{\beta}} + \frac{\{\epsilon_{\lambda} \beta - \lambda \epsilon_\beta\}(A-\pi)^2}{\beta \hat{\beta}} \\
    \E\left[\left\{\frac{\hat{\lambda}\hat{v}(O)}{\hat{\beta}} - \frac{\lambda v(O)}{\beta} \right\}^2 \Big| X\right] 
    &= \frac{\hat{\lambda}^2}{\hat{\beta}^2}\epsilon_\pi^2\{2\beta + \epsilon_\pi^2\} \\ 
    &+ \frac{\{\epsilon_{\lambda} \beta - \lambda \epsilon_\beta\}^2 \E\{(A-\pi)^4|X\}}{\beta^2\hat{\beta}^2} \\ 
    &+ \frac{2\hat{\lambda}\epsilon_\pi[2\E\{(A-\pi)^3|X\} + \epsilon_\pi \beta]\{\epsilon_{\lambda} \beta - \lambda \epsilon_\beta\}}{\beta\hat{\beta}^2}
\end{align*}
Note that Cauchy-Schwarz implies $\E\{(A-\pi)^3|X\}^2 \leq \E\{(A-\pi)^4|X\}\beta < K^2$, and $\lambda^2 \leq \beta \Var(Y|X) < K^2$.

% For the latter we write
% \begin{align*}
%     \hat{v}(O) - v(O) &= \frac{(A-\hat{\pi})^2}{\hat{\beta}} \left\{\epsilon_\lambda  - \epsilon_\beta\frac{\lambda}{\beta} \right\} -\epsilon_\pi\frac{\lambda}{\beta}\{2(A-\hat{\pi}) + \epsilon_\pi\}
% \end{align*}

% Choosing $\epsilon > 0$ such that $\hat{\beta} \geq \epsilon$ almost surely
% \begin{align*}
%     \E\left[\left\{\hat{g}(O) - g(O)\right\}^2\right] \leq \epsilon^{-2} \E\left[\left\{\hat{g}(O) - g(O)\right\}^2\right]
% \end{align*}

\subsection{Proof of Theorem \ref{asym_theorem_alse2}}

Let
\begin{align*}
\gamma &\equiv P\left[ \{Y-\mu(X)\}\{A-\pi(X)\}\right] \\
\hat{\gamma} &\equiv \mathbb{P}_n\left[ \{Y-\hat{\mu}(X)\}\{A-\hat{\pi}(X)\}\right] \\
\eta &\equiv P\left[ \{A-\pi(X)\}^2\right] \\
\hat{\eta} &\equiv \mathbb{P}_n\left[\{A-\hat{\pi}(X)\}^2\right] \\
\phi_{\gamma}(O) &\equiv \{Y-\mu(X)\}\{A-\pi(X)\} - \gamma \\
\hat{\phi}_{\gamma}(O) &\equiv \{Y-\hat{\mu}(X)\}\{A-\hat{\pi}(X)\} - \hat{\gamma} \\
\phi_{\eta}(O) &\equiv \{A-\pi(X)\}^2 - \eta \\
\hat{\phi}_{\eta}(O) &\equiv \{A-\hat{\pi}(X)\}^2 - \hat{\eta}
\end{align*}
Without making any restrictions we write
\begin{align*}
    \hat{\gamma} - \gamma &= (\mathbb{P}_n -P)\{\phi_{\gamma}(O)\} + R_n + H_n \\
    R_n &\equiv \hat{\gamma} - \gamma + P\{\hat{\phi}_{\gamma}(O)\} \\
    H_n &\equiv (\mathbb{P}_n -P)\{\hat{\phi}_{\gamma}(O) - \phi_{\gamma}(O)\}.
\end{align*}
We will show that the remainder therm $R_n=o_P(n^{-1/2})$ and the empirical process term $H_n=o_P(n^{-1/2})$. Since $P\{\phi_{\gamma}(O)\} = 0$, we therefore obtain the RAL results
\begin{align*}
    \hat{\gamma} - \gamma &= \mathbb{P}_n\{\phi_{\gamma}(O)\} + o_p(n^{-1/2}) \\
    \hat{\eta} - \eta &= \mathbb{P}_n\{\phi_{\eta}(O)\} + o_p(n^{-1/2})
\end{align*}
where the result for $\hat{\eta}$ follows as a special case ($Y=A$) of the result for $\hat{\gamma}$. Our goal is to consider the estimator $\hat{\Psi} = \hat{\gamma}/\hat{\eta}$, which estimates $\Psi = \gamma/\eta$. It follows by algebraic manipulations that,
\begin{align*}
\sqrt{n}(\hat{\Psi} - \Psi) &= \frac{\eta}{\hat{\eta}} \left[ \sqrt{n}\mathbb{P}_n \{\phi_{\Psi}(O)\} + o_p(1) \right]
\end{align*}
where $ \phi_{\Psi}(o)  = \{\phi_{\gamma}(o) - \Psi \phi_{\eta}(o)\}/\eta$ is the influence curve of $\Psi$. Next we use Slutsky's Theorem and the fact that $\hat{\eta}/\eta$ converges to 1 in probability, to write,
\begin{align*}
\lim_{n\to \infty} \sqrt{n}(\hat{\Psi} - \Psi) &= \lim_{n\to \infty} \sqrt{n}\mathbb{P}_n \{\phi_{\Psi}(O)\}
\end{align*}
which is the desired result.

\textbf{The remainder term}
To simplify notation, we will omit function arguments, e.g. $\mu = \mu(X)$ with similar for $\hat{\mu},\pi,\hat{\pi}$. Evaluating the remainder gives
\begin{align*}
R_n &= \E\left[ (Y-\hat{\mu})(A-\hat{\pi})-  (Y-\mu)(A-\pi) \right] \\
&= \E\left[ (\mu-\hat{\mu})(\pi-\hat{\pi})\right].
\end{align*}
By the Cauchy-Schwarz inequality $R_n^2 \leq \E\left[ (\mu-\hat{\mu})^2\right] \E\left[(\pi-\hat{\pi})^2\right]$, which is $o_P(n^{-1})$ under (A1).

\textbf{The empirical process term}

First write the empirical process term as the sum
\begin{align*}
    H_n &= (\mathbb{P}_n - P)\{\gamma - \hat{\gamma}\}\\
    &+ (\mathbb{P}_n - P)\left[\{Y-\mu(X)\}\{\hat{\pi}(X) - \pi(X)\}\right] \\
    &+ (\mathbb{P}_n - P)\left[\{\hat{\mu}(X) - \mu(X)\}\{A-\pi(X)\}\right] \\
    &+ (\mathbb{P}_n - P)\left[\{\hat{\mu}(X) - \mu(X)\}\{\hat{\pi}(X) - \pi(X)\}\right]
\end{align*}
Note that the first term is zero since $(\mathbb{P}_n - P)\{\gamma - \hat{\gamma}\} = (\gamma - \hat{\gamma})(\mathbb{P}_n - P)\{1\} = 0$. When the Donsker condition holds, then, by Lemma 19.24 of \cite{Vaart2013}, the remaining terms are $o_P(n^{-1/2})$ provided (i) that 
\begin{align*}
&\E\left[\Var(Y|X)\{\hat{\pi}(X) - \pi(X)\}^2\right] \\
&\E\left[\{\hat{\mu}(X) - \mu(X)\}^2\Var(A|X)\right] \\
&\E\left[\{\hat{\mu}(X) - \mu(X)\}^2\{\hat{\pi}(X) - \pi(X)\}^2\right]
\end{align*}
Similarly, under sample splitting then by Chebyshev's inequality, (i) is also sufficient conditions for $H_n$ to be $o_P(n^{-1/2})$. We claim that (i) holds since $\hat{\mu}$ and $ \hat{\pi}$ are consistent and $\Var(A|X) < K$ and $\Var(Y|X) < K$.

\section{Additional illustrated results}
\label{appendix:iwpc_extra}

\begin{table}[htb]
\caption{Least squares estimands applied to IWPC data, using the discrete Super Learner algorithm for model fitting. Results indicate the point estimates, its standard error, and 95\% Wald confidence interval, all in units of INR/(mg/week). P-values correspond to a Wald test of the null hypothesis that the estimand is zero.}
\centering
\begin{tabular}{|l|l|l|l|l|l|}
\hline
Estimand            & Algorithm                 & Estimate & SE & CI & p\\ \hline
$\Psi$     &noSS    & 1.87$\times10^{-3}$ & 6.32$\times10^{-4}$ &(0.633$\times10^{-3}$, 3.11$\times10^{-3}$) & 0.003       \\ \hline
$\Psi$     &SS    & 1.85$\times10^{-3}$ & 6.26$\times10^{-4}$ & (0.623$\times10^{-3}$, 3.08$\times10^{-3}$) & 0.003   \\ \hline
$\psi$     &noSS-A   & -0.105 & 0.956 & (-1.98, 1.77) & 0.91      \\ \hline
$\psi$     &SS-A   & -0.269 & 0.236 & (-0.732, 0.193) & 0.25   \\ \hline
$\psi$     &noSS-B   & 1.46$\times10^{-3}$ & 8.13$\times10^{-4}$ & (-0.133$\times10^{-3}$, 3.05$\times10^{-3}$) & 0.07       \\ \hline
$\psi$     &SS-B   & 1.37$\times10^{-3}$ & 8.43$\times10^{-4}$ & (-0.282$\times10^{-3}$, 3.02$\times10^{-3}$) & 0.10   \\ \hline
\end{tabular}
\label{IWPCresults2}
\end{table}

\begin{figure}[htb]
    \centering
    \includegraphics[width=\linewidth]{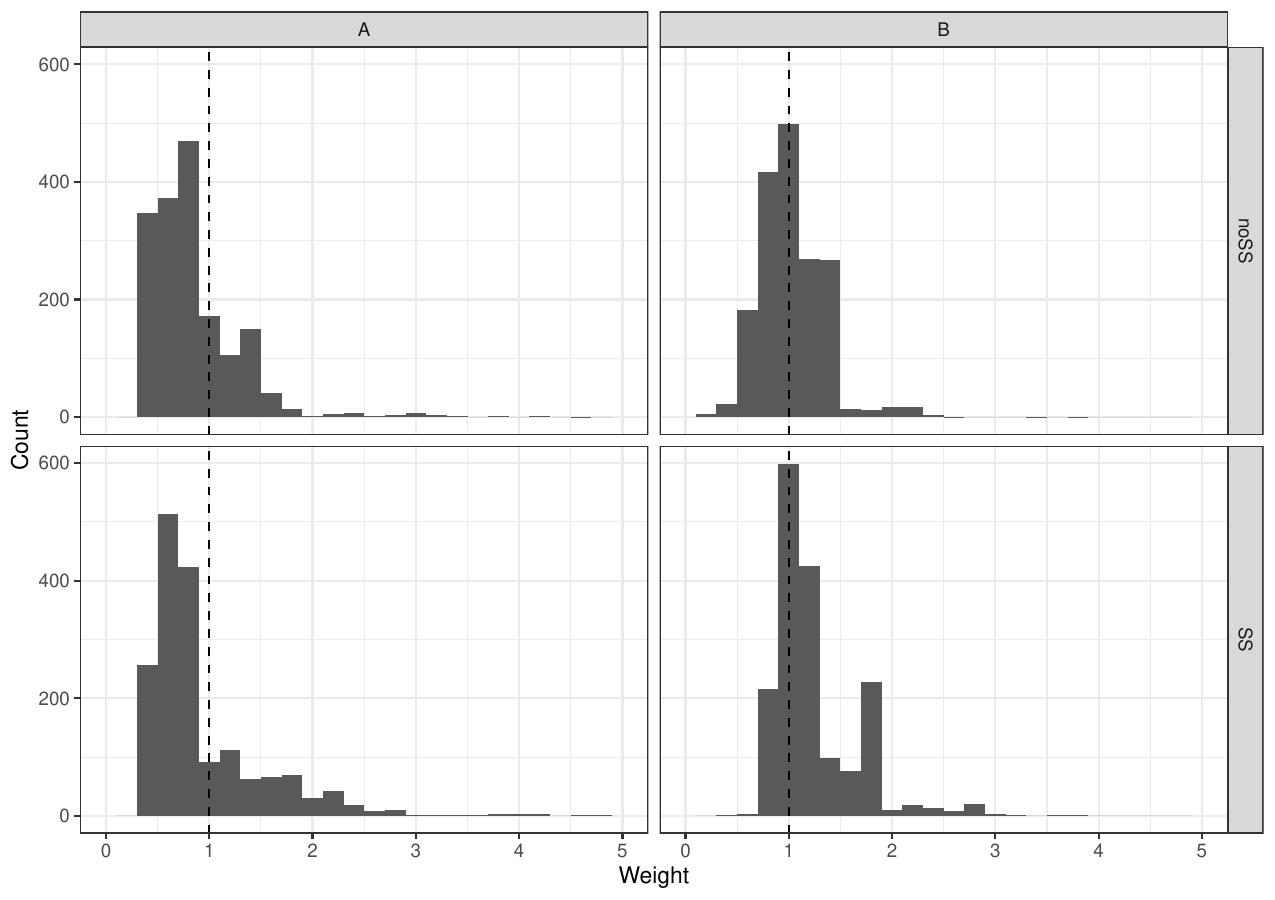}
    \caption{Least squares estimand weights approximated using the location-scale procedure described in Supplement \ref{weight_approximations}. This procedure uses estimates of the conditional mean and variance or $A$ given $X$, which are obtained using the algorithms in Section \ref{proposed_algos} and a discrete super learner for model fitting. Rows and columns refer to different algorithms as labeled.}
    \label{fig:discrete-superlearner-weights}
\end{figure}